\documentstyle[11pt]{article}
\newfam\msbfam
\font\tenmsb=msbm10    \textfont\msbfam=\tenmsb
\font\sevenmsb=msbm7   \scriptfont\msbfam=\sevenmsb
\font\fivemsb=msbm5    \scriptscriptfont\msbfam=\fivemsb
\def\Bbb{\fam\msbfam \tenmsb}

\def\rr{{\Bbb R}}
\def\rz{{{\rr}^n}}
\def\zz{{\Bbb Z}}
\def\nn{{\Bbb N}}

\def\fz{\infty}
\def\az{\alpha}
\def\supp{{\rm{\ supp\ }}}

\def\ez{\epsilon}
\def\bz{\beta}
\def\pz{\partial}
\def\tz{\theta}
\def\sz{\sigma}
\def\vz{\varphi}
\def\dz{\delta}
\def\gz{\gamma}

\def\lz{\lambda}
\def\oz{\Omega}
\def\supp{{\rm supp}}

\def\wt{\widetilde}
\def\wz{\omega}
\def\l{\left}
\def\r{\right}

\def\dsum{\displaystyle\sum}

\def\dint{\displaystyle\int}
\def\dfrac{\displaystyle\frac}
\def\dsup{\displaystyle\sup}

\def\dinf{\displaystyle\inf}

\newtheorem{thm}{\hskip\parindent Theorem}

\newtheorem{lem}{\hskip\parindent Lemma}

\newtheorem{prop}{\hskip\parindent Proposition}

\newtheorem{cor}{\hskip\parindent Corollary}

\newtheorem{defn}{\hskip\parindent Definition}

\begin{document}

\baselineskip=15pt
\renewcommand{\arraystretch}{2}
\arraycolsep=1.2pt

\title{Weighted norm inequalities for  Schr\"odinger type operators
 \footnotetext{ \hspace{-0.65
cm} 2000 Mathematics Subject  Classification: 42B25, 42B20.\\
The  research was supported  by the NNSF (10971002) of China.\\}}

\author{ Lin Tang }
\date{}
\maketitle

{\bf Abstract}\quad  Let $L=-\Delta+V$ be a Schr\"{o}dinger
operator, where $\Delta $ is the Laplacian operator on $\rz$,
while nonnegative potential $V$ belongs to the reverse H\"{o}lder
class. In this paper, we establish  the weighted norm inequalities
for some Schr\"odinger type operators, which include Riesz transforms and fractional integrals  and their commutators. These results  generalize substantially some well-known results.

\bigskip

\begin{center}{\bf 1. Introduction }\end{center}

In this paper, we consider the  Sch\"odinger differential operator
$$L=-\Delta+V(x)\ {\rm on}\ \rz,\ n\ge 3,$$where  $V(x)$
is a  nonnegative potential satisfying  certain reverse H\"older
class.

 We say a nonnegative locally $L^q$
integral function $V(x)$ on $\rr^n$ is said to belong to
$B_q(1<q\le \fz)$ if there exists $C>0$ such that the reverse
H\"older inequality
$$\l(\dfrac 1{|B(x,r)|}\dint_{B(x,r)} V^q(y)dy\r)^{1/q}\le C\l(\dfrac
1{|B(x,r)|}\dint_{B(x,r)} V(y)dy\r) \eqno(1.1)$$ holds for every
$x\in\rz$ and $0<r<\fz$, where $B(x,r)$ denotes the ball centered
at $x$ with radius $r$. In particular, if $V$ is a nonnegative
polynomial, then $V\in B_\fz$.  It is worth pointing out that the
$B_q$ class is that, if $V\in B_q$ for some $q>1$, then there
exists $\ez>0$, which depends only $n$ and the constant $C$ in
(1.1), such that $V\in B_{q+\ez}$. Throughout this paper, we
always assume that $0\not\equiv V\in B_{n/2}$.

The study of Schr\"odinger operator $L=-\triangle+V$ recently
attracted much attention; see \cite{bhs1,bhs2,dz,dz1,glp,s1,z}.  Shen \cite{s1} considered $L^p$
estimates for Schr\"odinger type operators $L$ with certain potentials
which include Schr\"odinger Riesz transforms $R_j^L=\frac{\pz}{\pz
x_j}L^{-\frac 12},\ j=1,\cdots,n$.
Very recently, Bongioanni, etc,  \cite{bhs1} proved
$L^p(\rz)(1<p<\fz)$ boundedness for commutators of Riesz transforms associated with Schr\"odinger operator with $BMO_\tz(\rho)$ functions which include the class $BMO$ function, and Bongioanni, etc,  \cite{bhs2} established the weighted boundedness for
Riesz transforms and fractional integrals associated with Schr\"odinger operator with weight $A_p^{\rho,\tz}$ class which includes the Muckenhoupt weight class. Naturally, it will be a very interesting problem to ask whether we can establish the weighted boundedness for commutators
of some Schr\"odinger type operators with $BMO_\tz(\rho)$ functions and weight $A_p^{\rho,\tz}$ class.

In this paper, we  give a confirm  answer. In order to answer question above, it seems that we can not adapt the methods from \cite{bhs1,bhs2}, so we need to use some new thoughts to overcome this obstacle in this paper. In fact, we establish a new Fefferman-Stein inequality and weighted inequalities for new maximal operators.  It is worth pointing out that our methods are more general than these in \cite{bhs1,bhs2}, since we can consider more general Schr\"odinger type operators by using our methods.

The paper is organized as follows. In Section 2,  we give some notation and basic results, these basic results play a crucial role in this paper.  In Section 3, we establish weighted norm inequalities for some Schr\"odinger type operators. In section 4,
we established the weighted boundedness for commutators of  Riesz transforms and fractional integrals associated with Schr\"odinger  operators.

Throughout this paper, we let $C$ denote   constants that are
independent of the main parameters involved but whose value may
differ from line to line. By $A\sim B$, we mean that there exists
a constant $C>1$ such that $1/C\le A/B\le C$.

\bigskip

\begin{center} {\bf 2. Some  notation and basic results}\end{center}

We first recall some notation.  Given $B=B(x,r)$
and $\lz>0$, we will write $\lz B$ for the $\lz$-dilate ball,
which is the ball with the same center $x$ and with radius $\lz
r$. Similarly, $Q(x,r)$ denotes the cube centered at $x$ with the
sidelength $r$ (here and below only cubes
with sides parallel to the coordinate axes are considered), and $\lz Q(x,r)=Q(x,\lz r)$.
 Given a Lebesgue
measurable set $E$ and a weight $\wz$, $|E|$ will denote the
Lebesgue measure of $E$ and $\wz(E)=\int_E\wz dx$.
$\|f\|_{L^p(\wz)}$ will denote $(\int_\rz |f(y)|^p\wz(y)dy)^{1/p}$
for $0< p<\fz$.

 The function $m_V(x)$ is defined by
$$\rho(x)=\dfrac 1{m_V(x)}=\dsup_{r>0}\l\{r:\ \dfrac
{1}{r^{n-2}}\dint_{B(x,r)}V(y)dy\le 1\r\}.$$Obviously,
$0<m_V(x)<\fz$ if $V\not=0$. In particular, $m_V(x)=1$ with $V=1$
and $m_V(x)\sim (1+|x|)$ with $V=|x|^2$.

\begin{lem}\label{l2.1.}\hspace{-0.1cm}{\rm\bf 2.1(\cite{s1}).}\quad
There exists $l_0>0$ and $C_0>1$such that $$\dfrac
1{C_0}\l(1+|x-y|m_V(x)\r)^{-l_0}\le \dfrac{m_V(x)}{m_V(y)}\le
C_0\l(1+|x-y|m_V(x)\r)^{l_0/(l_0+1)}.
$$ In particular, $m_V(x)\sim m_V(y)$ if $|x-y|<C/m_V(x)$.
\end{lem}

In this paper, we write $\Psi_\tz(B)=(1+r/\rho(x_0))^\tz$, where
 $\tz> 0$, $x_0$ and $r$ denotes the center and radius of $B$ respectively.

 A weight will always mean a positive function which is
locally integrable. As in \cite{bhs2}, we say that a weight $\wz$ belongs to the
class $A^{\rho,\tz}_p$ for $1<p<\fz$, if there is a constant $C$ such that
for all ball
  $B=B(x,r)$
$$\l(\dfrac 1{\Psi_\tz(B)|B|}\dint_B\wz(y)\,dy\r)
\l(\dfrac 1{\Psi_\tz(B)|B|}\dint_B\wz^{-\frac 1{p
-1}}(y)\,dy\r)^{p-1}\le C.$$
 We also
say that a  nonnegative function $\wz$ satisfies the $A^{\rho,\tz}_1$
condition if there exists a constant $C$ for all balls $B$
$$M_{V}^\tz(\wz)(x)\le C \wz(x), \ a.e.\ x\in\rz.$$
 where
$$M_{V}^\tz f(x)=\dsup_{x\in B}\dfrac 1{\Psi_\tz(B)|B|}\dint_B|f(y)|\,dy.$$
   Since $\Psi_\tz(B)\ge 1$, obviously,
$A_p\subset A_p^{\rho,\tz}$ for $1\le p<\fz$, where $A_p$ denote
the classical Muckenhoupt weights; see \cite{gr} and \cite{m}. We
will see that $A_p\subset\subset A_p^{\rho,\tz}$ for $1\le p<\fz$ in
some cases.  In fact, let $\tz>0$ and $0\le\gz\le\tz$, it is easy
to check that $\wz(x)=(1+|x|)^{-(n+\gz)}\not\in A_\fz=\bigcup_{p\ge 1}A_p$ and
$\wz(x)dx$ is not a doubling measure, but
$\wz(x)=(1+|x|)^{-(n+\gz)}\in A_1^{\rho,\tz}$ provided that  $V=1$ and
$\Psi_\tz(B(x_0,r))=(1+r)^\tz$.

We remark that  balls can be replaced by cubes  in definitions of
$A_p^{\rho,\tz}$ for $p\ge 1$ and $M_{V,\tz}$, since
$\Psi_\tz(B)\le  \Psi_\tz(2B)\le 2^\tz \Psi_\tz(B)$.
When $V=0$ and $\tz=0$, we
denote $M_{0,0}f(x)$ by $Mf(x)$( the standard Hardy-Littlewood
maximal function). It is easy to see that $|f(x)|\le M_{V}^\tz
f(x)\le Mf(x)$ for $a.e.\ x\in\rz$ and $\tz\ge 0$. For convenience, in the rest of this paper, fixed $\tz\ge 0$, we always assume that
$\Psi(B)$ denotes $\Psi_\tz(B)$ and $A_p^\rho$ denotes $A_p^{\rho,\tz}$.
\begin{lem}\label{l2.2.}\hspace{-0.1cm}{\rm\bf 2.2.}\quad Let $1<p<\fz$, then
\begin{enumerate}
\item[(i)]If $ 1\le p_1<p_2<\fz$, then $A_{p_1}^\rho\subset
A_{p_2}^\rho$. \item[(ii)] $\wz\in A_p^\rho$ if and only if
$\wz^{-\frac 1{p-1}}\in A_{p'}^\rho$, where $1/p+1/p'=1.$
\item[(iii)] If $\wz\in A_p^\rho$ for $1\le p<\fz$, then
$$\dfrac 1{\Psi(Q)|Q|}\dint_Q|f(y)|dy\le C\l(\dfrac
1{\wz(5Q)}\dint_Q|f|^p\wz(y)dy\r)^{1/p},$$ where
$\wz(E)=\int_E\wz(x)dx$. In particular, let $f=\chi_{E}$ for any
measurable set $E\subset Q$,
$$\dfrac {|E|}{\Psi(Q)|Q|}\le C\l(\dfrac
{\wz(E)}{\wz(5Q)}\r)^{1/p}.$$
\end{enumerate}
\end{lem}
{\it Proof.}\quad (i) and (ii) can be easily obtained by the definition of $A_p^\rho$. We only prove (iii).  In fact,
$$\begin{array}{cl}
\dfrac 1{\Psi(Q)|Q|}\dint_Q|f(y)|dy&=\dfrac
1{\Psi(Q)|Q|}\dint_Q|f(y)|\wz^{\frac 1p}(y)\wz^{-\frac 1p}(y)dy\\
&\le \l(\dfrac 1{\Psi(Q)|Q|}\dint_Q|f(y)|^p\wz(y)dy\r)^{\frac
1p}\l(\dfrac
1{\Psi(Q)|Q|}\dint_Q\wz^{-\frac 1{p-1}}(y)dy\r)^{\frac {p-1}p}\\
&\le C\l(\dfrac 1{\Psi(Q)|Q|}\dint_Q|f(y)|^p\wz(y)dy\r)^{\frac
1p}\\
&\qquad\qquad\times\l(\dfrac
1{\Psi(5Q)|5Q|}\dint_{5Q}\wz^{-\frac 1{p-1}}(y)dy\r)^{\frac {p-1}p}\\
&\le C\l(\dfrac 1{\Psi(Q)|Q|}\dint_Q|f(y)|^p\wz(y)dy\r)^{\frac
1p}\l(\dfrac
1{\Psi(5Q)|5Q|}\dint_{5Q}\wz(y)dy\r)^{-\frac 1p}\\
&\le C\l(\dfrac 1{\wz(5Q)}\dint_Q|f|^p\wz(y)dy\r)^{1/p}.
\end{array}$$
Thus, (iii) is proved.\hfill$\Box$

We also need the dyadic maximal operator $M^\triangle_{V,\eta} f(x)$ with
$0<\eta<\fz$  defined by
$$M_{V,\eta}^\triangle f(x):=\dsup_{x\in   Q(dyadic\ cube)}\dfrac
1{\Psi(Q)^\eta|Q|}\dint_Q|f(x)|\,dx.$$
Let $0<\eta<\fz$, the
dyadic sharp maximal operator $M^\sharp_{V,\eta}f(x)$ is defined
by
$$\begin{array}{cl}
M_{V,\eta}^\sharp f(x)&:=\dsup_{x\in Q,r<\rho(x_0)}\dfrac
1{|Q|}\dint_{Q(x_0,r)}|f(y)-f_Q|\,dy+ \dsup_{x\in Q,r\ge \rho(x_0)}\dfrac
1{\Psi(Q)^\eta|Q|}\dint_{Q(x_0,r)}|f|\,dy\\
&\simeq\dsup_{x\in Q,r< \rho(x_0)}\dinf_C \frac
1{|Q|}\int_{Q(x_0,r)}|f(y)-C|\,dy+ \dsup_{x\in Q,r\ge\rho(x_0)}\dfrac 1{\Psi(Q)^\eta |Q|}\dint_{Q(x_0,r)}|f|\,dx
\end{array}$$
where $Q_{x_0}'s$ denote dyadic cubes $Q(x_0,r)$ and $f_Q=\frac
1{|Q|}\int_Q f(x)dx$.

A variant of dyadic maximal operator and dyadic sharp maximal
operator $$M^\triangle_{\dz,\eta}
f(x)=M^\triangle_{V,\eta}(|f|^\dz)^{1/\dz}(x)$$ and
$$M_{\dz,\eta}^\sharp
f(x)=M^\sharp_{V,\eta}(|f|^\dz)^{1/\dz}(x),$$
which will become the main tool in our scheme.

 Bongioanni, etc, \cite{bhs2} proved the following Lemma.
\begin{lem}\label{l2.3.}\hspace{-0.1cm}{\rm\bf 2.3.}\quad
If $\wz\in
A_\fz^\rho=\bigcup_{p\ge 1}A_p^\rho$,  then there exists constants $\wt C>0$ and
$\dz_1>0$, such that for any $Q=Q(x_0,r)\subset\rz$ with $r<
\rho(x_0)$ such that for any measurable $E\subset Q$
$$\dfrac{\wz(E)}{\wz(Q)}\le \wt C\l(\dfrac{|E|}{|Q|}\r)^{\dz_1}.$$
\end{lem}
\begin{lem}\label{l2.4.}\hspace{-0.1cm}{\rm\bf 2.4.}\quad
Let $0<\eta<\fz$ and  $f$ be a  locally integrable function on
$\rz$, $\lz>0$, and $\Omega_\lz=\{x\in \rz:\
M^\triangle_{V,\eta}f(x)>\lz\}$. Then $\Omega_\lz$ may be written
as a disjoint union of dyadic cubes $\{Q_j\}$ with
\begin{enumerate}
\item[(i)] $\lz<(\Psi(Q_j)^\eta|Q_j|)^{-1}\dint_{Q_j}|f(x)|\,dx,$
\item[(ii)] $|Q_j|^{-1}\dint_{Q_j}|f(x)|\,dx\le
2^n(8nC_0)^{(l_0+2)\tz\eta}\lz,$ for each cube $Q_j=Q_j(x_j,r_j)$ with $r_j<\rho(x_j)$. This has the immediate consequences:
\item[(iii)] $|f(x)|\le \lz \ {\rm for}\ a.e\ x\in
\rz\setminus\bigcup_jQ_j$
 \item[(iv)]$|\Omega_\lz|\le
\lz^{-1}\dint_\rz|f(x)|\,dx.$
\end{enumerate}
\end{lem}
The proof follows from the same argument of Lemma 1 in page 150 in
\cite{s}.

By Lemmas 2.3 and 2.4, we  establish the following  weighted``good
$\lz$" inequality.
\begin{thm}\label{t2.1.}\hspace{-0.1cm}{\rm\bf 2.1.}\quad
Let $\wz\in A_\fz^\rho$ and $0<\eta<\fz$. For a locally integrable
function $f$, and for $b$ and $\gz$ positive
$\gz<b<b_0=(8nC_0)^{-(l_0+2)\tz\eta}$, we have the following
inequality
$$\wz(\{x\in \rz: M^\triangle_{V,\eta} f(x)>\lz, M^\sharp_{V,\eta}f(x)\le\gz
\lz\})\le \wt Ca^{\dz_1}\wz(\{x\in \rz:\  M^\triangle_{V,\eta}
f(x)>b\lz\})\eqno(2.1)$$ for all $\lz>0$, where $a=2^n\gz/(1-\frac
b{b_0})$, $C_0$ is defined in Lemma 2.1, $\wt C$ and $\dz_1$ are defined in Lemma 2.3.
\end{thm}
{\it Proof}.\quad  We may assume that the set $\{x:\
M^\triangle_{V,\eta} f(x)>b \lz\}$ has finite measure, otherwise
the inequality (2.1) is obvious. From Lemma 2.4, then this set is
the union of disjoint maximal cubes $\{Q_j\}$. We let $Q=Q(x_0,r)$
denote one of these cubes.  We consider two cases about sidelength
$r$, that is, $r< 1/m_V(x_0)$ and  $r\ge 1/m_V(x_0)$.

Case 1. When $r< 1/m_V(x_0)$, let $\wt Q\supset Q$ be the parent
of $Q$, by the maximality of $Q$ we have $|f|_{\wt Q}\le
b\lz\Psi(\wt Q)^\eta\le b\lz/b_0$ by Lemma 2.1. So far all $x\in Q$
for which $ M_{V,\eta}^\triangle f(x)>\lz$, it follows that
$M_{V,\eta}^\triangle (f\chi_ Q)](x)>\lz$, and also that
$M_{V,\eta}^\triangle [(f-f_{\wt Q})\chi_Q](x)>(1-b/b_0)\lz$. By
the weak type (1,1) of $M_{V,\eta}^\triangle$(see (iv) of Lemma
2.4), we have
$$\begin{array}{cl}
|\{x\in Q: M^\triangle_{V,\eta} f(x)>t,
M^\sharp_{V,\eta}f(x)\le\gz t\}|&\le \dfrac
1{(1-b/b_0)\lz}\dint_Q|f-f_{\wt
Q}|dx\\
&\le \dfrac 1{(1-b/b_0)t}\dint_{\wt Q}|f-f_{\wt
Q}|dx\\
&\le \dfrac {|\wt Q|}{(1-b/b_0)t}\dinf_{x\in
Q}M_{V,\eta}^\sharp f(x)\\
&\le \dfrac {2^n\gz| Q|}{1-b/b_0},
\end{array}$$
if the set in question is not empty. So
$$|\{x\in Q: M^\triangle_{V,\eta} f(x)>\lz, M^\sharp_{V,\eta}f(x)\le\gz
\lz\}|\le a|Q|.$$ From this  and by Lemma 2.3, in the case
$r<1/m_V(x_0)$, we have
$$\wz(\{x\in Q: M^\triangle_{V,\eta} f(x)>\lz, M^\sharp_{V,\eta}f(x)\le\gz
\lz\})\le \wt C a^{\dz_1}\wz(Q).\eqno(2.2)$$
 Case 2. When $r\ge 1/m_V(x_0)$, note that
 $$b\lz<\dfrac 1{\Psi(Q)^\eta|Q|}\dint_Q|f(y)|dy\le \dinf_{x\in Q}M_{V,\eta}^\sharp
 f(x)\le\gz \lz,$$
but $\gz<b$, hence, the set in question is  empty. Thus (2.2)
holds for any $Q$, and hence (2.1). The proof of Theorem 2.1 is
complete.\hfill$\Box$

As a consequence of Theorem 2.1, we have the following result.

\begin{cor}\label{c2.1.}\hspace{-0.1cm}{\rm\bf 2.1.}\quad
Let $0<p,\ \eta,\ \dz<\fz$ and $\wz\in A_\fz$. There exists a
positive constant $C$ such that
$$\dint_\rz M^\triangle_{\dz,\eta}
 f(x)^p\wz(x)dx\le C \dint_\rz M^\sharp_{\dz,\eta}
f(x)^p\wz(x)dx.$$ Let  $\vz: (0,\fz)\to(0,\fz)$ be a doubling
function. Then there exists a positive constant $C$ such that
$$\sup_{\lz>0}\vz(\lz)\wz(\{x\in\rz:\ M^\triangle_{\dz,\eta}f(x)
>\lz\})\le C  \sup_{\lz>0}\vz(\lz)\wz(\{x\in\rz:\ M^\sharp_{\dz,\eta}f(x)
>\lz\})$$ for any smooth function $f$
for which the left handside is finite.
\end{cor}
To establish weighted inequality for fractional integrals, we need introduce $A_{(p,q)}^\rho$. We say that a
weight $\wz$ belongs to the class $A_{(p,q)}^\rho$ for $1\le p< \fz$
and $1\le q<\fz$. Let $p'=p/(p-1)$, if there is a constant $C$
such that for any cube
  $Q=Q(x,r)$
$$\l(\dfrac 1{\Psi(Q)|Q|}\dint_Q[\wz(y)]^q\,dy\r)^{1/q}
\l(\dfrac 1{\Psi(Q)|Q|}\dint_Q[\wz(y)]^{-p'}\,dy\r)^{1/p'}\le C.$$
Obviously, $\wz^{1/p}\in A_{(p,p)}^\rho\Leftrightarrow \wz\in
A_{p}^\rho$ for $1\le p<\fz$.

Next, let $0\le\bz<n$, we give a result about the operator
$M_{\bz,\wz}$ defined by
$$M_{\bz,\wz}(f)(x)=\dsup_{x\in B}\dfrac
1{\wz(5B)^{1-\bz/n}}\dint_B|f(x)|\wz(x)dx.$$ In the rest of this
section, we write $M_{\wz}f(x)=M_{0,\wz}f(x).$
\begin{lem}\label{l2.5.}\hspace{-0.1cm}{\rm\bf 2.5.}\quad
Let $0\le\bz<n, 1\le p<n/\bz$ and $1/q=1/p-\bz/n$. If $\wz\in
A_\fz^\rho$ , then
$$\wz(\{x\in\rz:\ M_{\bz,\wz}f(x)>\lz\})\le C\l(\dfrac
{\|f\|_{L^p(\wz)}}\lz\r)^{q},\ \forall\lz>0,\ \forall f\in
L^p(\wz).$$ In particular, from (3.1) and using Marcinkiewicz
interpolation theorem, then for $1< p<n/\bz$ and $1/q=1/p-\bz/n$
so that
$$\|M_{\bz,\wz}f\|_{L^q(\wz)}\le C\|f\|_{L^p(\wz)}.$$
\end{lem}
 {\it Proof}.\quad  We  set $x\in
E_\lz=\{x\in\rz:\ M_{\bz,\wz}f(x)>\lz\}$ with any $\lz>0$, then,
there exists a ball $B_x\ni x$ such that
$$\dfrac
1{\wz(5B_x)^{1-\bz/n}}\dint_{B_x}|f(y)|\wz(y)dy>\lz.\eqno(2.3)$$
Thus, $\{B_x\}_{x\in E_\lz}$ covers $E_\lz$. By Vitali lemma,
there exists a class disjoint cubes $\{B_{xj}\} $ such that
$\bigcup B_{xj}\subset E_\lz\subset \bigcup 5B_{xj}$ and
$$\wz(E_\lz)\le \dsum_j\wz(5B_{xj}).\eqno(2.4)$$
From (2.3), we have
$$\lz<\dfrac
1{\wz(5B_x)^{1/q}}\l(\dint_{B_x}|f(y)|^p\wz(y)dy\r)^{1/p}.$$ From
this and by (2.4), note that $p/q\le 1$, we get
$$\begin{array}{cl} \wz(E_\lz)^{p/q}&\le
\dsum_j\wz(5B_{xj})^{p/q}\le \dfrac C{\lz^p}\dsum_j\dint_{B_{xj}}|f(y)|^p\wz(y)dy\\
&=\dfrac C{\lz^p}\dint_{\bigcup_jB_{xj}}|f(y)|^p\wz(y)dy\le \dfrac
C{\lz^p}\dint_\rz|f(y)|^p\wz(y)dy.\end{array}$$ Thus, Lemma 2.5 is
proved.\hfill$\Box$

The fractional maximal operator $M_{\bz,V}$ is defined by
$$M_{\bz,V}f(x)=\dsup_{x\in Q}\dfrac 1{(\Psi(Q)|Q|)^{1-\frac \bz n}}\dint_Q|f(y)|\,dy,$$
where $0\le \bz<n$, we let $M_{V}$ denote $M_{0,V}$.

From (iii) of Lemma 2.2, we know that for $1\le p<\fz$ and $\wz\in A_p^\rho$
$$M_Vf(x)\le C(M_{\wz}(|f|^p)(x))^{1/p},\quad x\in\rz.$$
From this and using Lemma 2.5, we can get the following result.
\begin{prop}\label{p2.1.}\hspace{-0.1cm}{\rm\bf 2.1.}\quad
Let $1<p <\fz$ and suppose that $\wz\in A_p^\rho$. If $p<p_1<\fz$,
then the equality
$$\dint_\rz|M_Vf(x)|^{p_1}\wz(x)dx\le
C\dint_\rz|f(x)|^{p_1}\wz(x)dx.$$ Further, let $1\le p<\fz$,
$\wz\in A_p^\rho$   if and only if
$$\wz(\{x\in\rz:\ M_{V}f(x)>\lz\})\le \dfrac
{C}{\lz^p}\dint_\rz|f(x)|^p\wz(x)dx.$$
\end{prop}
For the fractional maximal operator $M_{\bz,V}$, we have:
\begin{prop}\label{p2.2.}\hspace{-0.1cm}{\rm\bf 2.2.}\quad
Let $0<\bz<n, 1\le p<n/\bz$ and $1/q=1/p-\bz/n$. If $\wz\in
A_{(p,q)}^\rho$, then
$$\l(\dint_{\{x\in\rz:\ M_{\bz,V}f(x)>\lz\})}[\wz(y)]^qdy\r)^{1/q}\le \dfrac
{C}{\lz}\l(\dint_\rz|f(x)|^p\wz(x)dx\r)^{1/p}.$$
\end{prop}
{\it Proof}.\quad Fix $M>0$ and let $E_{\lz, M}=\{|x|<M:
M_{\bz,V}f(x)>\lz\}$. For each $x\in E_{\lz,M}$ there is a cube
$Q\ni x$ such that
$$(\Psi(Q)|Q|)^{\frac \bz n-1}\dint_Q|f(y)|dy>\lz.\eqno(2.5)$$
Pick a sequence $\{Q_k\}$ of these cubes such that
$E_{\lz,M}\subset\bigcup Q_k$ and no point of $\rz$ is in more
than $L$ of these cubes where $L$ depends only on $n$ (see
\cite{mw}). Note that $p/q<1$, by (2.5), we then have
$$\begin{array}{cl}
\l(\dint_{E_{\lz, M}}[\wz(x)]^qdx\r)^{p/q}&\le
\l(\dsum_k\dint_{Q_k}[\wz(x)]^qdx\r)^{p/q}\le
\dsum_k\l(\dint_{Q_k}[\wz(x)]^qdx\r)^{p/q}\\
&\le C \dsum_k\l(\dint_{Q_k}[\wz(x)]^qdx\r)^{p/q}\l(
\lz^{-1}(\Psi(Q_k)|Q_k|)^{\frac \bz n-1}\dint_{Q_k}|f(x)|dx\r)^p\\
&\le C
\dsum_k\l(\dint_{Q_k}[\wz(x)]^qdx\r)^{p/q}\lz^{-p}(\Psi(Q_k)|Q_k|)^{1-p-p/q}\\
&\qquad\qquad\qquad\times\l(
\dint_{Q_k}|f(x)\wz(x)|^pdy\r)\l(\dint_{Q_k}[\wz(x)]^{-p'}dx\r)^{p/p'}\\
&\le C\lz^{-p}\dsum_k \dint_{Q_k}|f(x)\wz(x)|^pdx\\
&\le C\lz^{-p} \dint_\rz|f(x)\wz(x)|^pdx.
\end{array}$$
Using the monotone convergence theorem, we can obtain the desired
result.\hfill$\Box$

 Next we will  establish the weighted strong type $(p,q)$ for a variant maximal operator $M_{V,\eta}$ for
$0<\eta<\fz$ as follows
$$M_{\bz,V,\eta}f(x)=\dsup_{x\in Q}\dfrac 1{(\Psi(Q))^{\eta}(\Psi(Q)|Q|)^{1-\frac \bz n}}\dint_Q|f(y)|\,dy.$$
\begin{thm}\label{t2.2.}\hspace{-0.1cm}{\rm\bf 2.2.}\quad
Let $0\le\bz<n, 1< p<n/\bz, 1/p+1/p'=1$ and $1/q=1/p-\bz/n$. If $\wz\in
A_{(p,q)}^\rho$ and $\eta\ge(1-\bz/n)p'/q$, then there exists a
constant $C>0$ such that
$$\|M_{\bz,V,\eta}f\|_{L^q(\wz^q)}\le C\|f\|_{L^p(\wz^p)}.$$
\end{thm}
{\it Proof}.\quad  We shall adapt an argument in \cite{j1}. Note that $\wz\in A_{(p,q)}^\rho\Leftrightarrow \wz^q\in
A_{1+q/p'}^\rho$, where  $p'=p/(p-1)$. Let $\gz=1+q/p',\
\gz'=\gz/(\gz-1)$ and $\nu=\wz^q$. Set $\sz=\nu^{-\frac
1{\gz-1}}$, then $\sz\in A_{\gz'}^\rho$. Let $f\in L^p(\wz^p)$, for
any $k\in\zz$, taking any compact set $K_k\subset\{x\in\rz:
2^k<M_{\bz,V,\eta}f(x)\le 2^{k+1}\}$, for any $x\in K_k$, $\exists
Q_x\ni x$ such that
$$ 2^{k+1}\ge\dfrac 1{(\Psi(Q_x))^{\eta}(\Psi(Q_x)|Q_x|)^{1-\bz/n}}\dint_{Q_x}f(y)|dy>2^k.$$
We take finite cover $\{Q_j^k\}$ from the cover $\{Q_x\}_{x\in
K_k}$ of $K_k$. Write
$$E_1^k=Q_1^k\bigcap K_k,\ E_j^k=(Q_j^k-\bigcup_{i<j}Q_j^k)\bigcap
K_k,\ j>1,$$ then $\{E_j^k\}$ is a disjoint collection on $j$ for
fixed $k$, and $K_k=\bigcup_j E_j^k$,  we then get
$$\begin{array}{cl}
\dint_{\bigcup_kK_k}|M_{\bz,V,\eta}f|^q\wz^q dx&\le
C\dsum_{k,j}2^{kq}\nu(E_j^k)\\
&\le C\dsum_{k,j}\nu(E_j^k)\l(\dfrac 1{(\Psi(Q_j^k))^{\eta}(\Psi(Q_j^k)|Q_j^k|)^{1-\bz/n}}\dint_{Q_j^k}|f(y)|dy\r)^q\\
\end{array}$$
$$\begin{array}{cl}
&= C\dsum_{k,j}\nu(E_j^k)\l(\dfrac
{\sz(5Q_j^k)^{1-\bz/n}}{(\Psi(Q_j^k))^{\eta}(\Psi(Q_j^k)|Q_j^k|)^{1-\bz/n}}\r)^q\\
&\qquad\qquad\times \l(\dfrac
1{\sz(5Q_j^k)^{1-\bz/n}}\dint_{Q_j^k}|f(y)|\sz^{-1}\sz dy\r)^q.
\end{array}\eqno(2.6)$$
Define the measure $\mu$ on $\zz\times\zz_+$ by
$$\mu:\ (k,j)\to\mu_{k,j}=\nu(E_j^k)\l(\dfrac
{\sz(5Q_j^k)^{1-\bz/n}}{(\Psi(Q_j^k))^{\eta}(\Psi(Q_j^k)|Q_j^k|)^{1-\bz/n}}\r)^q
.$$ Write
$$\Gamma(\lz)=\l\{(k,j):\ \l(\dfrac
1{\sz(5Q_j^k)^{1-\bz/n}}\dint_{Q_j^k}|f(y)|\sz^{-1}\sz
dy\r)^q>\lz\r\},\ \lz>0,$$
$$G(\lz)=\bigcup\{Q_j^k:\ (k,j)\in\Gamma(\lz)\}.$$
Since $\nu\in A_\gz^\rho$, then
$$\dsup_Q\l(\dfrac{\sz(5Q)}{\Psi(Q)|Q|}\r)^\gz\l(\dfrac
{\nu(5Q)}{\Psi(Q)|Q|}\r)^{\gz'}\le  C.$$ Note that
$\gz=(1-\bz/n)q$ and $\eta\ge(1-\bz/n)p'/q$, we then  obtain
$$\begin{array}{cl}
\mu_{k,j}&=\nu(E_j)\l(\dfrac
{\sz(5Q_j^k)^{1-\bz/n}}{(\Psi(Q_j^k))^{\eta}(\Psi(Q_j^k)|Q_j^k|)^{1-\bz/n}}\r)^q\le C\nu(E_j^k)\l(\dfrac {|Q_j^k|}{\nu(5Q_j^k)}\r)^{\gz'}\\
 & \le
C\nu(E_j^k)\l(\dfrac 1{\nu(5Q_j^k)}\dint_{Q_k^j}\nu^{-1}\nu dy\r)^{\gz'}\\
&\le C\nu(E_j^k)\dinf_{x\in
Q_j^k}M_{\nu}(\nu^{-1}\chi_{Q_j^k})^{\gz'}(x)\\
&\le
C\dint_{E_j^k}M_{\nu}(\nu^{-1}\chi_{Q_j^k})^{\gz'}(x)\nu(x)dx.
\end{array}$$
Since $M_{\nu}$ is bounded on $L^{\gz'}(\nu)$ by Lemma 2.5, we
have
$$\begin{array}{cl}
\mu(\Gamma (\lz))&=\dsum_{(k,j)\in\Gamma (\lz)}\mu_{k,j}\le C
\dsum_{(k,j)\in\Gamma
(\lz)}\dint_{E_j^k}M_{\nu}(\nu^{-1}\chi_{G(\lz)})^{\gz'}(x)\nu(x)dx\\
&\le \dint_{G(\lz}M_{\nu}(\nu^{-1}\chi_{G(\lz)})^{\gz'}(x)\nu(x)dx\\
&\le C\dint_{G(\lz)}\nu^{1-\gz'}dx=\sz(G(\lz)).
\end{array}$$
Note that
$$G(\lz)\subset\{x\in\rz:\ (M_{\bz,\sz}(f\sz^{-1}))^q(x)>\lz\},$$ so
$$\sz(G(\lz))\le\sz(\{x\in\rz:\ (M_{\bz,\sz}(f\sz^{-1}))^q(x)>\lz\}).$$
Hence, by Lemma 2.5, the right side of (2.6) equal to
$$\begin{array}{cl}
\dint_0^\fz\mu(\Gamma(\lz))d\lz&\le C\dint_0^\fz\sz(G(\lz)d\lz\\
&\le C\dint_0^\fz\sz(\{x\in\rz:\ (M_{\bz,\sz}(f\sz^{-1}))^q(x)>\lz\})d\lz\\
&=C\dint_\rz (M_{\bz,\sz}(f\sz^{-1}))^q\sz dx\\
&\le C\l(\dint_\rz
|f|^p\sz^{1-p}dx\r)^{q/p}\\
&=C\l(\dint_\rz|f\wz|^p dx\r)^{q/p}.
\end{array}$$
Thus,
$$\l(\dint_\rz(M_{\bz,V,\eta}f)^q\wz^q dx\r)^{1/q}\le C\l(\dint_\rz |f|^p\wz^p dx\r)^{1/p}.$$
 Theorem 2.2 is proved.\hfill$\Box$

We next recall some basic definitions and facts about Orlicz
spaces, referring to \cite{r} for a complete account.

 A function $B(t): [0,\fz)\to [0,\fz)$ is called a Young function if
it is continuous, convex, increasing and satisfies $\Phi(0)=0$ and
$B\to \fz$ as $t\to\fz$. If $B$ is a Young function, we define the
$B$-average of a function $f$ over a cube $Q$ by means of the
following Luxemberg norm:
$$\|f\|_{B,Q}=\dinf\l\{\lz>0:\ \dfrac
1{|Q|}\dint_QB\l(\dfrac {|f(y)|}\lz\r)\,dy\le 1\r\}.$$ The
generalized H\"older's inequality
$$\dfrac
1{|Q|}\dint_B|fg|\,dy\le \|f\|_{B,Q}\|g\|_{\bar B,Q}$$ holds,
where $\bar B$ is the complementary Young function associated to
$B$. And we define the corresponding maximal function
$$M_B f(x)=\dsup_{Q:x\in Q}\|f\|_{B,Q}$$
and for $0<\eta<\fz$ and $0\le\bz<n$
$$M_{B,\bz,V,\eta} f(x)=\dsup_{Q:x\in Q}\Psi(Q)^{-\eta}(\Psi(Q)|Q|)^{\bz/n}\|f\|_{B,Q}.$$
In particular, if $\bz=0$, we denote $M_{B,0,V,\eta} f(x)$ by
$M_{B,V,\eta} f(x)$.

 The in example that we are going to use is
$B(t)=t(1+log^+t)$ with the maximal function denoted by
$M_{LlogL}$. The complementary Young function is given by $\bar
B(t)\approx e^t$ with the corresponding maximal function denoted
by $M_{exp L}$.

\begin{center} {\bf 3. Weighted norm inequalities for Schr\"odinger type operators }\end{center}
 We first consider a class Schr\"odinger type operators such as
$\nabla(-\Delta+V)^{-1}\nabla$, $\nabla(-\Delta+V)^{-1/2}$,
$(-\Delta+V)^{-1/2}\nabla$  with $V\in B_n$, $(-\Delta+V)^{i\gz}$  with $\gz\in\rr$ and $V\in B_{n/2}$, and $\nabla^2(-\Delta+V)^{-1}$ with $V$ is a nonnegative polynomial, are standard Calder\'on-Zygmund operators,  in particular,  the kernels $K$ of operators above all satisfy
the following conditions for some $\dz_0>0$ and  any $l\in \nn_0=\nn\bigcup\{0\}$,
$$|K(x,y)|\le\dfrac {C_l}{(1+|x-y|(m_V(x)+m_V(y)))^l}\dfrac 1{|x-y|^n}\eqno(3.1)$$
and
$$|K(x+h,y)-k(x,y)|+|K(x,y+h)-K(x,y)|\le\dfrac {C_l}{(1+|x-y|(m_V(x)+m_V(y)))^l}\dfrac {|h|^{\dz_0}}{|x-y|^{n+\dz_0}},\eqno(3.2)$$
whenever $x,y,h\in\rz$, and $|h|<|x-y|/2$, and $m_V(x)$ is defined in Section 2.

\begin{thm}\label{t3.1.}\hspace{-0.1cm}{\rm\bf 3.1.}\quad
Let $T$ denote the operators above. Let $1<p <\fz$ and suppose that $\wz\in A_p^\rho$. Then
$$\dint_\rz|Tf(x)|^p\wz(x)dx\le
C\dint_\rz|f(x)|^p\wz(x)dx.$$ Further, suppose that $\wz\in
A_1^\rho$. Then, there exists a constant $C$ such that for all
$\lz>0$
$$\wz(\{x\in\rz:\ Tf(x)>\lz\})\le \dfrac
C\lz\dint_\rz|f(x)|\wz(x)dx.$$
\end{thm}
We remark that the weighted  boundedness of $\nabla(-\Delta+V)^{-1/2}$,
$(-\Delta+V)^{-1/2}\nabla$  with $V\in B_n$ is proved in \cite{bhs2}.
 We can prove Theorem 3.1
by using the similar proof of   Theorems 3.4 below. We omit the details here. 

Next we give a result of maximal Schr\"odinger type  operators.

\begin{thm}\label{t3.2.}\hspace{-0.1cm}{\rm\bf 3.2.}\quad
Let $0<p, \eta <\fz$ and suppose that $\wz\in A_\fz^\rho$. Then
$$\dint_\rz|T^*f(x)|^p\wz(x)dx\le
C_{p,\eta}\dint_\rz|M_{V,\eta}f(x)|^p\wz(x)dx$$ and
$$\dsup_{\lz>0}\lz\wz(\{x\in\rz:\ T^*f(x)>\lz\})\le C\dsup_{\lz>0}\lz\wz(\{x\in\rz:\ M_{V,\eta}f(x)>\lz\}),$$
 where  the
maximal operator $T^*$ defined by
$$T^*f(x):=\dsup_{\ez>0}|T_\ez f(x)|=\dsup_{\ez>0}\l|\dint_{|y-x|>\ez}K(x,y)f(y)dy\r|.$$
\end{thm}
To prove theorem 3.2, we need the following lemma.
\begin{lem}\label{l3.1.}\hspace{-0.1cm}{\rm\bf 3.1.}\quad
For any a  ball $B=B(x_0,r)$, if $r\ge 1/m_V(x_0)$, then the ball
$B$ can be decomposed into finite disjoint cubes $\{Q_i\}_{i=1,m}$
such that $B\subset \bigcup_i^{m} Q_i\subset 2\sqrt{n}B$ and
$r_i/2\le 1/m_V(x)\le 2\sqrt{n}C_0r_i$ for some $x\in
Q_i=Q(x_i,r_i)$, where $C_0$ is  same as Lemma 2.1.
\end{lem}
{\it Proof}.\quad In fact, let Q be a cube with center at $x_0$ and
sidelength $2r$. Obviously, $B\subset Q\subset 2\sqrt{n}B$. If
there exist a point $x\in Q$ such that $2r/2\le 1/m_V(x)$, by
Lemma 2.1 and $r\ge 1/m_V(x_0)$, we then know that $1/m_V(x)\le
2\sqrt{n}C_02r$, thus $Q$ will be satisfied. Otherwise, we split
$Q$ into $2^n$ disjoint cubes $O_i$ with sidelength $r$. If  $Q_i$
satisfies $r/2\le 1/m_V(x)\le 2\sqrt{n}C_0r$ for some $x\in Q_i$,
we keep it, otherwise, we continuous split $Q_i$ as above. From
Lemma 2.1, we know that
$1/m(x)>(1+2r\sqrt{n}/m_V(x_0))^{-l_0}/(C_0m_V(x_0))$ for all
$x\in Q$. Therefore, the splitting steps must be finite. Thus, we
can obtain finite disjoint cubes $Q_i$ such that
$Q=\bigcup_{i=1}^{m}Q_i$ and $r_i/2\le 1/m_V(x)\le 2C_0r_i$ for
some $x\in Q_i=Q(x_i,r_i)$. Obviously, these cubes $Q_i's$ are
just what we need. Thus, Lemma 3.1 is proved.\hfill$\Box$

\medskip

{\it Proof of Theorem 3.2}.\quad The set $\oz=\{x\in\rz: \ T^*f(x)>t\}$ is open.
Therefore we can decompose it as a disjoint union $\oz=\bigcup
Q_j$ of Whitney cubes: they are mutually disjoint and $2\
diam(Q_j)\le dist(Q_j,\oz^c)\le 8\ diam(Q_j)$. Moreover, the
family $4Q_j$ is almost disjoint with constant $4^n$, and
obviously $4Q_j\subset\oz$.

To prove Theorem 3.2, we only need to show that
$$\wz(\{x\in \rz: T^*f(x)>(1+\bz)t\ {\rm and}\ M_{V,\eta}f(x)\le\gz
t\})\le a\wz(\{x\in \rz:\ T^*f(x)>t\})\eqno(3.3)$$ for all $t>0$
and $a<(1+\bz)^{-1-\frac 1p}$.

Let $Q_j=Q(x_j,r_j)$, we set $E_1=\{j: r_j\le 1/m_V(x_j)\}$ and
$E_2=\{j: r_j>1/m_V(x_j)\}$. From the proof of Lemma 3.1, we know
that for any $j\in E_2$, the cube $Q_j$ can be decomposed into
finite disjoint cubes $\{Q^i_j\}_{i=1}^{j_m}$ such that $Q_j=
\bigcup_{i=1}^{j_m} Q^i_j$ and $r_{ji}/2\le 1/m_V(x)\le
2\sqrt{n}C_0r_{ji}$ for some $x\in Q^i_j$. We are going to show
that, given $\bz>0$ and $0<\az<1$, there exists
$\gz=\gz(\bz,\az,n)$ such that for all $j\in E_1$
$$|\{x\in Q_j: T^*f(x)>(1+\bz)t\ {\rm and}\ M_{V,\eta}f(x)\le\gz
t\}|\le 4^n\az|Q_j|$$  and $j\in E_2$ and $Q_j=\bigcup_{i=1}^{j_m}
Q^i_j$ such that
$$|\{x\in Q_j^i: T^*f(x)>(1+\bz)t\ {\rm and}\ M_{V,\eta}f(x)\le\gz
t\}|\le 4^n\az|Q_j^i|$$ From these and using Lemma 3.1 with
$q_0=(20nC_0)^{(l_0+2)\az}$, we have
$$\wz(\{x\in Q_j: T^*f(x)>(1+\bz)t\ {\rm and}\ M_{V,\eta}f(x)\le\gz
t\})\le C\az^{\dz_1}\wz(Q_j),\ j\in E_1\eqno(3.4)$$  and $j\in
E_2$ and $Q_j=\bigcup_{i=1}^{j_m} Q^i_j$ such that
$$\wz(\{x\in Q_j^i: T^*f(x)>(1+\bz)t\ {\rm and}\ M_{V,\eta}f(x)\le\gz
t\})\le C\az^{\dz_1}\wz(Q_j^i).\eqno(3.5)$$ Summing over $j$ and
$i$, we get
$$\wz(\{x\in \rz: T^*f(x)>(1+\bz)t\ {\rm and}\ M_{V,\eta}f(x)\le\gz
t\})\le C\az^{\dz_1}\wz(\oz).$$ Choose $\az$ so that
$C\az^{\dz_1}<(1+\bz)^{-1-\frac 1p}$, we would finally obtain
(3.5).

It remains to prove (3.4) and (3.5). We first consider the term
(3.4). Fix $j$ and set $Q=Q_j=Q(x_0,r)$. Assume that there exists
$\bar x\in Q$  so that $M_{V,\eta}f(\bar x)\le\gz t$(if not, the
set appearing in (3.4) would be empty). Let $z\in \oz^c$ such that
$dist(z,Q)=dist(Q,\oz^c)$.  Note that
$$Q\subset P=Q(\bar x,\frac 52r)\subset 4Q\subset Q_z=Q(z,18r).$$
Set $f_1=f\chi_{Q_z}$ and $f_2=f-f_1$.  Note that $r<1/m_V(x_0)$
implies  $\Psi(Q_z)\sim 1$,  for $x\in Q$, we then have
$$\begin{array}{cl}
|T_\ez f_1(x)|&\le|T_\ez (f\chi_P)(x)|+\dfrac
C{r^n}\dint_{Q_z}|f(y)|dy\\
&\le |T^* (f\chi_P)(x)|+CM_{V,\eta}f(\bar x)\\
&\le |T^* (f\chi_P)(x)|+C\gz t.
\end{array}$$
Hence,
$$|T_\ez f(x)|\le|T_\ez f_2(x)|+|T^* (f\chi_P)(x)|+C\gz t.$$
By direction computation, we obtain
$$|T_\ez f_2(x)-T_\ez f_2(z)|\le C M_{V,\eta}f(\bar x)$$ and
$$|T_\ez f_2(z)|\le T^*f(z)\le t.$$
Then
$$T^*f(x)\le T^*(f\chi_P)(x)+(1+C\gz)t,\quad x\in Q.$$
Take $\gz$ so that $2C\gz\le \bz$, we then have
$$\{x\in Q_j: T^*f(x)>(1+\bz)t\ {\rm and}\ M_{V,\eta}f(x)\le\gz
t\}\subset\{x\in Q: T^*(f\chi_P)(x)>\frac\bz 2 t\}.$$By the weak
type (1,1) of   $T^*$, we get
$$\begin{array}{cl}
|\{x\in Q: T^*(f\chi_P)(x)>\frac\bz 2 t\}|&\le\dfrac C{\bz
t}\dint_P|f(y)|dy\\
&\le \dfrac{C|Q|}{\bz t}\dfrac 1{|4Q|}\dint_{4Q}|f(y)|dy\\
&\le \dfrac{C|Q|}{\bz t}M_{V,\eta}f(\bar x)\\
&\le \dfrac{C\gz|Q|}{\bz}\le \az|Q|,
\end{array}$$
if $\gz$ is chosen small enough so that $C\bz^{-1}\gz\le\az$.

Finally, we consider the term (3.5). Similar to (4.6), fix $j,i$
and set $Q=Q_j^i$ and $r=l(Q)$. Assume that there exists $\bar
x\in Q$  so that $M_{V,\eta}f(\bar x)\le\gz t$. Note that
$$Q\subset P=Q(\bar x,\frac 52r)\subset 4Q\subset Q_{\bar x}=Q(\bar x,18r).$$
Set $f_1=f\chi_{Q_{\bar x}}$ and $f_2=f-f_1$. Then for $x\in Q$,
we have
$$\begin{array}{cl}
|T_\ez f_1(x)|&\le|T_\ez (f\chi_P)(x)|+\dfrac
C{r^n}\dint_{Q_{\bar x}}|f(y)|dy\\
&\le |T^* (f\chi_P)(x)|+CM_{V,\eta}f(\bar x)\\
&\le |T^* (f\chi_P)(x)|+C\gz t,
\end{array}$$
so
$$|T_\ez f(x)|\le|T_\ez f_2(x)|+|T^* (f\chi_P)(x)|+C\gz t.$$
By direction computation, we obtain
$$|T_\ez f_2(x)|\le C M_{V,\eta}f(\bar x).$$
Then
$$T^*f(x)\le T^*(f\chi_P)(x)+C\gz t,\quad x\in Q.$$
The rest proof is similar to that of (3.4), we omit the details.
Thus,  Theorem 3.2 is proved.\hfill$\Box$

Next we consider another class $V\in B_q$ for $n/2\le q$ for Riesz transforms associated to Schr\"odinger operators. Let $T_1=(-\triangle+V)^{-1}V,\
T_2= (-\triangle+V)^{-1/2} V^{1/2}$ and $T_3=(-\triangle+V)^{-1/2} \nabla$.

\begin{thm}\label{t3.3.}\hspace{-0.1cm}{\rm\bf 3.3.}\quad
Suppose $V\in B_q$ and $q\ge n/2$.  Then
\begin{enumerate}
\item[(i)] If $q'\le p<\fz$ and $\wz\in A_{p/q'}^\rho$,
$$\| T_1f\|_{L^p(\wz)}\le C\|f\|_{L^p(\wz)};$$
\item[(ii)] If $(2q)'\le p<\fz$ and $\wz\in A_{p/(2q)'}^\rho$,
$$\| T_2f\|_{L^p(\wz)}\le C\|f\|_{L^p(\wz)};$$
\item[(iii)] If $p_0'\le p<\fz$ and $\wz\in A_{p/p_0'}^\rho$, where $1/p_0=1/q-1/n$ and $n/2\le q<n$,
$$\| T_3f\|_{L^p(\wz)}\le C\|f\|_{L^p(\wz)}.$$
\end{enumerate}
\end{thm}
Let $T_1^*=V(-\triangle+V)^{-1}, T^*_2=  V^{1/2}(-\triangle+V)^{-1/2}$ and $T^*_3=\nabla(-\triangle+V)^{-1/2} $. By duality we can easily get the following results.
\begin{cor}\label{c3.1.}\hspace{-0.1cm}{\rm\bf 3.1.}\quad
Suppose $V\in B_q$ and $q\ge n/2$.  Then
\begin{enumerate}
\item[(i)] If $1< p\le q$ and $\wz^{-\frac 1{p-1}}\in A_{p'/q'}^\rho$,
$$\| T^*_1f\|_{L^p(\wz)}\le C\|f\|_{L^p(\wz)};$$
\item[(ii)] If $1< p\le 2q$ and $\wz^{-\frac 1{p-1}}\in A_{p'/(2q)'}^\rho$,
$$\| T^*_2f\|_{L^p(\wz)}\le C\|f\|_{L^p(\wz)};$$
\item[(iii)] If $1< p\le p_0$ and $\wz^{-\frac 1{p-1}}\in A_{p'/p'_0}^\rho$, where $1/p_0=1/q-1/n$ and $n/2\le q<n$,
$$\| T^*_3f\|_{L^p(\wz)}\le C\|f\|_{L^p(\wz)}.$$
\end{enumerate}
\end{cor}
We remark that the weighted $L^p(\wz)$ boundedness of $T_3,\ T^*_3$ is proved in \cite{bhs2}.

To prove Theorem 3.3, we need the following result.
\begin{lem}\label{l3.2.}\hspace{-0.1cm}{\rm\bf 3.2.}\quad
Let $q, p_0$ be same as in Theorem 3.3. There exist constants $C_N>0$, $\dz>0$ and $s_1,s_2, s_3$ such that  $s_1>q, s_2>2q, s_3>p_0$, such that, for  $\forall\ r>0,
 x, x_0\in \rz$ with $|x-x_0|\le r$, then
$$\dsum_{k=1}^\fz (2^kr)^{\frac n{s_i'}}\l(\dint_{2^kr\le|y-x_0|<2^{k+1}r}|K(x,y)-K(x_0,y)|^{s_i} dy\r)^{1/s_i}\le \dsum_{k=1}^\fz \dfrac {C_N }{2^{k\dz}(1+
m_V(x_0)2^kr)^N},$$ where $K_i$ denotes the kernels of $T_i$ defined as above, $i=1,2,3$, and $1/s'_i+1/s_i=1$ for $i=1,2,3$.
\end{lem}Lemma 3.2 is essentially proved in \cite{glp}.

\bigskip

{\it Proof of Theorem 3.3}.\quad  For convenience, Let $T$ denote these operators $T_1,T_2,T_3$,  $K,s$ denote the kernel $K_i$ of $T_i$ and  $s_i$ for $i=1,2,3$ respectively. From Corollary 2.1 and Proposition 2.1, and note that for any  $\eta>0$, $|f(x)|\le M^\triangle_{V}f(x),\ a.e.\ x\in\rz$, we need only
to show that
$$M_V^\sharp(Tf)(x)\le CM_{V}^{1/s}(|f|^s)(x), \ {\rm a.e\ }\ x\in\rz,\eqno(3.6)$$
holds for any $f(x)\in C_0^\fz(\rz)$.

We fix $x\in\rz$ and let $x\in Q=Q(x_0,r)$(dyadic cube).
Decompose $f=f_1+f_2$, where $f_1=f\chi_{\bar Q}$, where $\bar
Q=Q(x,8\sqrt{n}r)$. Let  $C_Q$ a constant to be
fixed along the proof.

To prove (3.3), we consider two cases about $r$, that is, $r<
\rho(x_0)$ and $r\ge \rho(x_0)$.

Case 1. when  $r< \rho(x_0)$.
  We then have
$$\begin{array}{cl}
\dfrac 1{|Q|}\dint_Q|Tf(y)-C_Q|\,dy
&\le\dfrac
1{|Q|}\dint_Q|T(f_1)(y)|\,dy\\
&\qquad+\dfrac
1{|Q|}\dint_Q|T(f_2)(y)-C_Q|\,dy\\
&=I+II.
\end{array}$$
 To deal with $I$,  note that $m_V(x)\sim m_V(x_0)$ for any $x\in \bar Q$ and
$\Psi(\bar Q)\sim 1$, by $L^s(\rz)$ boundedness of $T$ (see \cite{s1}), we then have
$$\begin{array}{cl}
I &\le C \l(\dfrac
1{|Q|}\dint_Q|Tf(y)|^{s}\,dy\r)^{1/s}\\
&\le C \l(\dfrac
1{|Q|}\dint_{\bar Q}|f(y)|^{s}\,dy\r)^{1/s}\le C M_{V}^{1/s}(|f|^s)(x),
\end{array}\eqno(3.7)$$ where $1/s'+1/s=1$.

Finally, for II we first  fix the value of $C_Q$ by taking
$C_Q=(T(f_2))_Q$, the average of $T(f_2)$ on $Q$. Let $Q_k=Q(x_0,2^{k+1}r)$. By Lemmas
2.1 and 3.2, we then have
$$\begin{array}{cl}
II&\le \dfrac C{|Q|^2}\dint_Q\dint_Q\dint_{\rz\setminus \bar
Q}|K(y,\wz)-K(z,\wz)||f(\wz)|d\wz dzdy\\
&\le \dfrac C{|Q|^2}\dint_Q\dint_Q\dint_{|x_0-\wz|>16r
}|K(y,\wz)-K(z,\wz)||f(\wz)|d\wz dzdy\\
&\le \dfrac
C{|Q|^2}\dint_Q\dint_Q\dsum_{k=2}^\fz\dint_{2^kr\le|x_0-\wz|<2^{k+1}r
}|K(y,\wz)-K(z,\wz)||f(\wz)|d\wz dzdy\\
&\le C_N \dsum_{k=2}^\fz \dfrac {C_N }{2^{k\dz}(1+
m_V(x_0)2^kr)^N}
\l(\dint_{Q_k}|f(\wz)|^{s}d\wz\r)^{1/s}\le C_N M_{V}(|f|^s)^{1/s}(x),
\end{array}\eqno(3.8)$$  where $N=(l_0+1)\tz\eta$.

Case 2. when  $r> \rho(x_0)$.
  We then have
$$\begin{array}{cl}
\dfrac 1{\Psi(Q)|Q|}\dint_Q|Tf(y)|\,dy
&\le\dfrac
1{\Psi(Q)|Q|}\dint_Q|T(f_1)(y)|\,dy+\dfrac
1{\Psi(Q)|Q|}\dint_Q|T(f_2)(y)|\,dy\\
&:=I_1+II_1.
\end{array}$$
 To deal with $I_1$,  by $L^s(\rz)$ boundedness of $T$ again, we then have
$$\begin{array}{cl}
I &\le C \l(\dfrac
1{\Psi(Q)^{s}|Q|}\dint_Q|Tf(y)|^{s}\,dy\r)^{1/s}\\
&\le C \l(\dfrac
1{\Psi(Q)^s|Q|}\dint_{\bar Q}|f(y)|^{s}\,dy\r)^{1/s}\le C M_{V}^{1/s}(|f|^s)(x).
\end{array}\eqno(3.9)$$
For II, by Lemma
2.1, we then have
$$\begin{array}{cl}
II&\le \dfrac C{|Q|}\dint_Q\dint_{\rz\setminus \bar
Q}|k(y,\wz)||f(\wz)|d\wz dy\\
&\le \dfrac C{|Q|}\dint_Q\dint_{|x_0-\wz|>16r
}|k(y,\wz)||f(\wz)|d\wz dy\\
&\le \dfrac
C{|Q|}\dint_Q\dsum_{k=2}^\fz\dint_{2^kr\le|x_0-\wz|<2^{k+1}r
}|k(y,\wz)||f(\wz)|d\wz dy\\
&\le C_N \dsum_{k=2}^\fz \dfrac {C_N }{(1+
m_V(x_0)2^kr)^N|Q_k|}
\l(\dint_{Q_k}|f(\wz)|d\wz\r)\le C_N M_{V}(|f|^s)^{1/s}(x),
\end{array}\eqno(3.10)$$  where $N=(l_0+1)\tz\eta+1$.

 From  (3.7)--(3.10), we get (3.6). Hence the proof is
 finished.\hfill$\Box$

Finally, we establish some weighted inequalities for fractional
integrals associated with Schr\"odinger operators
 defined by
$${\cal I}_\bz f(x) ={ L}^{-\bz/2}f(x)=
\dint_0^\fz e^{-t{ L}}f(x)t^{\bz/2-1}dt=\dint_\rz
k_\bz(x,y)f(y)dy\quad {\rm for}\quad 0<\bz<n.$$
Using Proposition 2.4 in \cite{dz1}, we can get the following
result for the fractional integral associated with Schr\"odinger
operator.
\begin{lem}\label{l3.3.}\hspace{-0.1cm}{\rm\bf 3.3.}\quad
If $V\in B_q(\rz), q\ge n/2$ and $0<\bz<n$, $k_\bz$ denotes the
kernel of the fractional integral $I_\bz$ as above, then there
exists $\dz_0=\dz_0(q)>0$  such that for every $l>0$ there is a
constant $C_l$ so that
$$|k_\bz(x,y)|\le\dfrac {C_l}{(1+|x-y|(m_V(x)+m_V(y)))^l}\dfrac 1{|x-y|^{n-\bz}}$$ and
$$|k_\bz(x+h,y)-k_\bz(x,y)|\le\dfrac {C_l}
{(1+|x-y|(m_V(x)+m_V(y))^l}\dfrac
{|h|^{\dz_0}}{|x-y|^{n-\bz+\dz_0}},$$ whenever $x,y,h\in\rz$ and
$|h|<|x-y|/2.$
\end{lem}

\begin{thm}\label{t3.4.}\hspace{-0.1cm}{\rm\bf 3.4.}\quad
Let $0<\bz<n, 1< p<\bz/n$ and $1/q=1/p-\bz/n$. If $\wz\in
A_{(p,q)}^\rho$, then
$$\l(\dint_\rz|I_\bz f(x)|^q\wz(x)^qdx\r)^{1/q}\le
C\l(\dint_\rz|f(x)|^p\wz(x)^pdx\r)^{1/p}.$$ Further, suppose that
$\mu=\wz^q\in A_1^\rho$ with $q=n/(n-\bz)$. Then, there exists a
constant $C$ such that for all $\lz>0$
$$\mu(\{x\in\rz:\ |I_\bz f(x)|>\lz\})^{1/q}\le \dfrac
C\lz\dint_\rz|f(x)|\wz(x)dx.$$
\end{thm}
Theorem 3.4 is proved in \cite{bhs2}. Here we give another proof.  We will see that Theorem 3.4
follows from Proposition 2.2, Theorem 2.2 and  Theorem 3.5 below. It is worth pointing out that Theorem 3.5 has own interesting.

\begin{thm}\label{t3.5.}\hspace{-0.1cm}{\rm\bf 3.5.}\quad
Let $0<q, \eta <\fz$ and suppose that $\wz\in A_\fz^\rho$, then
$$\dint_\rz|I_\bz f(x)|^q\wz(x)dx\le C\dint_\rz|M_{\bz,V,\eta}f(x)|^q\wz(x)dx$$ and
$$\dsup_{\lz>0}\lz^q\wz(\{x\in\rz:\ |I_\bz f(x)|>\lz\})\le C\dsup_{\lz>0}\lz^q\wz(\{x\in\rz:\ M_{\bz,V,\eta}f(x)>\lz\}),$$
\end{thm}
{\it Proof}.\quad To prove Theorem 3.3, from Corollary 2.1 and Theorem
2.2, we need only to show that for any $0<\eta<\fz$ and
$0<\dz<\eta/(\eta+1)$ such that
$$M^{\sharp}_{\dz,\eta}(I_\bz f)(x)\le C M_{\bz,V,\eta}(f)(x)), \ {\rm a.e\ }\ x\in\rz,\eqno(3.11)$$
Fix $x\in\rz$ and let $x\in Q=Q(x_0,r)$(dyadic cube). Decompose
$f=f_1+f_2$, where $f_1=f\chi_{\bar Q}$, where $\bar Q=Q(x,8\sqrt{n}r)$.

To prove (3.11), we consider two cases about $r$, that is, $r<
1/m_V(x_0)$ and $r\ge 1/m_V(x_0)$.

Case 1. when  $r< \rho(x_0)$.
  Let $C_Q=|(I_\bz f_2)_Q|$.
Since $0<\dz<1$, we then have
$$\begin{array}{cl}
\l(\dfrac 1{|Q|}\dint_Q||I_\bz f(y)|^\dz-C_Q^\dz|\,dy\r)^{1/\dz}&
\le \l(\dfrac
1{|Q|}\dint_Q||I_\bz f(y)|-|(I_\bz f_2)_Q||^\dz\,dy\r)^{1/\dz}\\
&\le \l(\dfrac
1{|Q|}\dint_Q|I_\bz f(y)-(I_\bz f_2)_Q|^\dz\,dy\r)^{1/\dz}\\
&\le C\l(\dfrac
1{|Q|}\dint_Q|I_\bz f_1(y)|^\dz\,dy\r)^{1/\dz}\\
&\qquad+C\l(\dfrac
1{|Q|}\dint_Q|I_\bz f_2(y)-(I_\bz f_2)_Q|^\dz\,dy\r)^{1/\dz}\\
&=I+II.
\end{array}$$
For I, we recall that $I_\bz$ is weak type $(1,n/(n-\bz))$. Note
that $m_V(x)\sim m_V(x_0)$ for any $x\in \bar Q$ and $\Psi(\bar
Q)\sim 1$, by Kolmogorov's inequality(see\cite{p}), we then have

$$\begin{array}{cl}
I&\le \dfrac C{|Q|^{1-\bz/n}}\|I_\bz( f_1)\|_{L^{\frac n{n-\bz},\fz}}\\
&\le \dfrac C{|\bar Q|^{1-\bz/n}}\dint_{\bar Q}|f(y)|\,dy\le CM_{\bz,V,\eta}f(x).
\end{array}\eqno(3.12)$$
For II, let $Q_k=Q(x_0,2^{k+1}r)$ and $\az=\eta+1$. By
Lemma 3.3, we then have
$$\begin{array}{cl}
II&\le \dfrac C{|Q|}\dint_Q|I_\bz(f_2)(y)-(I_\bz(f_2))_Q|\,dy\\
&\le \dfrac C{|Q|^2}\dint_Q\dint_Q\dint_{\rz\setminus \bar
Q}|k_\bz(y,\wz)-k_\bz(z,\wz)||f(\wz)|d\wz dzdy\\
&\le \dfrac C{|Q|^2}\dint_Q\dint_Q\dint_{|x_0-\wz|>2r
}|k_\bz(y,\wz)-k_\bz(z,\wz)|f(\wz)|d\wz dzdy\\
&\le \dfrac
C{|Q|^2}\dint_Q\dint_Q\dsum_{k=4}^\fz\dint_{2^kr\le|x_0-\wz|<2^{k+1}r
}|k_\bz(y,\wz)-k_\bz(z,\wz)||f(\wz)|d\wz dzdy\\
&\le C_l\dsum_{k=2}^\fz\dfrac {2^{-
k}}{(1+2^krm_V(x_0))^l(2^{k+1}r)^n}\dint_{Q_k}|f(\wz)|d\wz\\
&\le C_l\dsum_{k=1}^\fz\dfrac
{2^{-k}(1+2^krm_V(Q_k))^{\az\tz}}{(1+2^krm_V(x_0))^l}\dfrac
1{(1+2^krm_V(Q_k))^{\az\tz}|Q_k|^{1-\bz/n}}\dint_{Q_k}
|f(\wz)|d\wz\\
 &\le C_l\dsum_{k=1}^\fz\dfrac
{2^{-k}(1+2^krm_V(x_0))^{\az(l_0+1)\tz}}{(1+2^krm_V(x_0))^l}\dfrac
1{(1+2^krm_V(Q_k))^{\az\tz}|Q_k|^{1-\bz/n}}\dint_{Q_k}
|f(\wz)|d\wz\\
&\le C_l\dsum_{k=1}^\fz  2^{- k}M_{\bz,V,\eta}(f)(x)\le
C_lM_{\bz,V,\eta}(f)(x),
\end{array}\eqno(3.13)$$ if taking $l=\tz(l_0+1)\az$.

Case 2. When  $r\ge \rho(x_0)$, note that
$\az_1:=\eta/\dz\ge\eta+1$, we have
 $$\begin{array}{cl} \dfrac
C{\Psi(Q)^{\az_1}} \l(\dfrac 1{|Q|}\dint_Q|I_\bz
f(y)|^\dz\,dy\r)^{1/\dz}& \le \dfrac C{\Psi(Q)^{\az_1}}\l(\dfrac
1{|Q|}\dint_Q|I_\bz f_1(y)|^\dz\,dy\r)^{1/\dz}\\
&\qquad+\dfrac C{\Psi(Q)^{\az_1}}\l(\dfrac 1{|Q|}\dint_Q|I_\bz
f_2(y)|^\dz\,dy\r)^{1/\dz}\\
 &:=I_1+II_1.
\end{array}$$
For $I_1$, similar to $I$, we have
$$\begin{array}{cl}
I_1&\le \dfrac C{\Psi(Q)^{\az_1}}\dfrac 1{|Q|^{1-\bz/n}}\|I_\bz f_1\|_{L^{\frac n{n-\bz},\fz}}\\
&\le \dfrac C{\Psi(Q)^\eta(\Psi(Q)|\bar Q|)^{1-\bz/n}}\dint_{\bar Q}|f(y)|\,d\mu(y)\\
&\le CM_{\bz,V,\eta}f(x).
\end{array}\eqno(3.8)$$
Finally, for $II_1$,  by Lemma
3.3, we then have
$$\begin{array}{cl}
II_1&\le \dfrac C{|Q|}\dint_Q|I_\bz(f_2)(y)|\,dy\\
&\le \dfrac C{|Q|}\dint_Q\dint_{\rz\setminus \bar
Q}|k_\bz(y,\wz)||f(\wz)|d\wz dy\\
&\le \dfrac C{|Q|}\dint_Q\dint_{|x_0-\wz|>2r
}|k_\bz(y,\wz)|f(\wz)|d\wz dy\\
&\le \dfrac
C{|Q|}\dint_Q\dsum_{k=2}^\fz\dint_{2^kr\le|x_0-\wz|<2^{k+1}r
}|k_\bz(y,\wz)|f(\wz)|d\wz dy\\
&\le C_l\dsum_{k=2}^\fz\dfrac 1{(1+2^krm_V(x_0))^l|Q_k|^{1-\bz/n}}\dint_{Q_k}|f(\wz)|d\wz\\
&\le C_l\dsum_{k=1}^\fz\dfrac
{(1+2^krm_V(Q_k))^{\tz\az_1}}{(1+2^krm_V(x_0))^l}\dfrac
1{(1+2^krm_V(Q_k))^{\tz\az_1}|Q_k|^{1-\bz/n}}\dint_{Q_k}
|f(\wz)|d\wz\\
 &\le C_l\dsum_{k=1}^\fz\dfrac
{(1+2^krm_V(x_0))^{\tz(l_0+1)\az_1}}{(1+2^krm_V(x_0))^l}\dfrac
1{(1+2^krm_V(Q_k))^{\tz\az_1}|Q_k|^{1-\bz/n}}\dint_{Q_k}
|f(\wz)|d\wz\\
&\le C_l\dsum_{k=1}^\fz  2^{- k}M_{\bz,V,\eta}(f)(x)\le
C_lM_{\bz,V,\eta}(f)(x),
\end{array}\eqno(3.14)$$  if taking $l=\az_1(l_0+1)\tz+1$.

 From  (3.12)--(3.14), we get (3.11). Hence the proof is
 finished.\hfill$\Box$

\begin{center} {\bf 4. Commutators for Shr\"odinger type operators }\end{center}
Bongioanni, etc, \cite{bhs1} introduce a new space $BMO_\tz(\rho)$ defined by
$$\|f\|_{BMO_\tz(\rho)}=\dsup_{B\subset
\rz}\dfrac 1{\Psi(B)|B|}\dint_B|f(x)-f_B|dx<\fz,$$ where
$f_B=\frac 1{|B|}\int_Bf(y)dy$ and $\Psi(B)=(1+r/\rho(x_0))^\tz$, $B= B(x_0,r)$ and $\tz> 0$.

Bongioanni, etc, \cite{bhs1} proved the following result for $BMO_\tz(\rho)$.
\begin{prop}\label{p4.1.}\hspace{-0.1cm}{\rm\bf 4.1.}\quad
Let $\tz>0$ and $1\le s<\fz$. If $b\in BMO_\tz(\rho)$, then
$$\l(\dfrac 1{|B|}\dint_B|b-b_B|^s\r)^{1/s}\le c C_0^\tz s\|b\|_{BMO_\tz(\rho)}\l(1+\dfrac r{\rho(x)}\r)^{\tz'},$$ for all $B=B(x,r)$, with
$x\in\rz$ and $r>0$, where $\tz'=(l_0+1)\tz$ and $C_0$ is defined in Lemma 2.1 and $c$ is a constant depending only on $n$.
\end{prop}

Obviously, the classical $BMO$ is properly contained in $BMO_\tz(\rho)$; more examples see \cite{bhs1}. For convenience, we let $BMO(\rho)$ denote $BMO_\tz(\rho)$.
\begin{prop}\label{p4.2.}\hspace{-0.1cm}{\rm\bf 4.2.}\quad
If $f\in BMO(\rho)$, then there exist  positive constant $c_1$
and $c_2$ such that for every ball $B$  and every $\lz>0$, we have
$$|\{x\in B: |f(x)-f_B|>\lz\}|\le c_1|B|\exp\l\{-\frac
{c_2\lz}{\|f\|_{BMO(\rho)}\Psi_{\tz'}(B)}\r\},$$
where
$f_B=\frac 1{|B|}\int_Bf(y)dy$ and $\Psi_{\tz'}(B)=(1+r/\rho(x_0))^{\tz'}$, $B= B(x_0,r)$ and $\tz'=(l_0+1)\tz$.
\end{prop}
{\it Proof}.\quad We adapt the same argument of pages 145-146 in \cite{st}. We first assume $\|f\|_{BMO(\rho)}$ $\Psi_{\tz'}(B)=1$. We apply Chebysheff's inequality and Proposition 4.1, we obtain
$$ |\{x\in B: |f(x)-f_B|>\lz\}|\le (c C_0^\tz s)^s\lz^{-s}|B|$$ for $0<\lz<\fz,\ 1\le s<\fz$.

If $\lz\ge 2c C_0^\tz$, we take $s=\lz/(2c C_0^\tz)\ge 1$. Then
$$ |\{x\in B: |f(x)-f_B|>\lz\}|\le (1/2)^s|B|=e^{-c_1\lz}|B|$$
where $c_1=(2c C_0^\tz)^{-1}\ln 2$. However, if $\lz\le 2c C_0^\tz$, then $e^{-c_1\lz}\ge e^{-c_12(c C_0^\tz)}=1/2$, and
$$ |\{x\in B: |f(x)-f_B|>\lz\}|\le 2e^{-c_1\lz}|B|$$ in that range of $\lz$. Altogether then, if we drop the normalization on $f$ by replacing $f$ by
$f/(\|f\|_{BMO(\rho)}\Psi_{\tz'}(B))$, we can obtain the conclusion by taking $c_1=2$ and $c_2=(2c C_0^\tz)^{-1}\ln 2$.\hfill$\Box$

As a consequence of Proposition 4.2, we can then obtain the following
result, which is equivalent to Proposition 4.1.

\begin{prop}\label{p4.3.}\hspace{-0.1cm}{\rm\bf 4.3.}\quad
 Suppose that $f$ is in $BMO(\rho)$. There exist
positive constants $\gz$ and $C$ such that
$$\dsup_B\dfrac 1{|B|}\dint_B\exp\l\{\frac
\gz{\|f\|_{BMO(\rho)}\Psi_{\tz'}(B)}|f(x)-f_B|\r\}\,dx\le C.$$
\end{prop}
{\it Proof}. We choose $\gz= c_2/2$, where $ c_2$ is a constant
in Proposition 4.2. We then have
$$\begin{array}{cl}
\dint_B&\exp\l\{\dfrac
\gz{\|f\|_{BMO(\rho)}\Psi_{\tz'}(B)}|f(x)-f_B|\r\}\,dx\\
&=\dint_0^\fz |\{x\in B: \exp\l\{\frac
\gz{\|f\|_{BMO(\rho)}\Psi_{\tz'}(B)}|f(x)-f_B|\r\}>t\}|\,dt\\
&\le |B|+\dint_1^\fz |\l\{x\in B: |f(x)-f_B|>\dfrac
{\log {t}\|f\|_{BMO(\rho)}\Psi_{\tz'}(B)}\gz\r\}|\,dt\\
&\le |B|+ c_1|B|\dint_1^\fz \exp\l\{-\dfrac {c_2\log
t}\gz\r\}\,dt\\
&\le |B|+ C|B|\dint_1^\fz t^{-2}\,dt\\
&\le C|B|.
\end{array}$$
Thus, the proof  is complete.\hfill$\Box$

We first consider   commutators of fractional integrals associated with Sch\"odinger operators.
\begin{thm}\label{t4.1.}\hspace{-0.1cm}{\rm\bf 4.1.}\quad
Let $b\in BMO(\rho)$, $0<\bz<n, 1< p<n/\bz$ and $1/q=1/p-\bz/n$. If
$\wz\in A_{(p,q)}^\rho$, then there exists a constant $C$ such that
$$\l(\dint_\rz|[I_\bz,b] f(x)|^q\wz(x)^qdx\r)^{1/q}
\le C \|b\|_{BMO(\rho)}\l(\dint_\rz|f(x)|^p\wz(x)^pdx\r)^{1/p}.$$
\end{thm}
The weighted weak-type endpoint estimate for the commutator is the
 following.
\begin{thm}\label{t4.2.}\hspace{-0.1cm}{\rm\bf 4.2.}\quad
Let $b\in BMO(\rho)$, $0<\bz<n, 1< p<n/\bz$ and $1/q=1/p-\bz/n$. Let
$B(t)=t\log(e+t)$, $\Lambda(t)=[t\log(e+t^{\bz/n})]^{n/(n-\bz)}$
and $\Theta(t)=t^{1-\bz/n}\log(e+t^{-\bz/n})$. If $\wz\in A_1^\rho$,
then for any $\lz>0$
$$\wz(\{x\in\rz: |[I_\bz,b] f(x)|>\lz\})\le
C\Lambda\l(\dint_\rz
B\l(\dfrac{\|b\|_{BMO(\rho)}|f(x)|}\lz\r)\Theta(\wz(x))dx\r).$$
\end{thm}
To prove Theorems 4.1 and 4.2, we need a few of Lemmas which can
be of independent interest.
\begin{lem}\label{l4.1.}\hspace{-0.1cm}{\rm\bf 4.1.}\quad
Let $0\le\bz<n,\ 0<\eta<\fz$ and $M_{V,\eta/2}f$ be locally
integral. Then there exist positive constants $C_1$ and $C_2$
independent of $f$ and $x$ such that
$$C_1M_{\bz,V,\eta}M_{V,\eta+1}f(x)\le M_{L\log L,\bz, V,\eta+1}f(x)\le
C_2M_{\bz,V,\eta/2} M_{V,\eta/2}f(x).$$
\end{lem}
{\it Proof}.\quad We first prove $M_{L\log L,\bz, V,\eta+1}f(x)\le
CM_{\bz,V,\eta/2} M_{V,\eta/2}f(x)$. It suffices to show that for
any cube $Q\ni x$, there is a constant $C>0$ such that
$$\Psi(Q)^{-\eta-1}(\Psi(Q)|Q|)^{\bz/n}\|f\|_{L\log L,Q}\le C
\Psi(Q)^{-\eta/2}(\Psi(Q)|Q|)^{\bz/n-1}\dint_QM_{V,\eta/2}f(y)dy.$$
That is,
$$\|f\|_{L\log L,Q}\le C\Psi(Q)^{\eta/2}{|Q|}^{-1}\dint_QM_{V,\eta/2}f(y)dy.\eqno(4.3)$$
In fact, by homogeneity we can take $f$ with $\|f\|_{L\log L,Q}=1$
which implies
$$\begin{array}{cl}
1&\le \dfrac C{|Q|}\dint_Q|f(y)|(1+\log^+(|f(y)|))dy\\
&\le \dfrac C{|Q|}\dint_Q|f(y)|\dint_1^{|f(y)|+e}\dfrac {dt}tdy\\
&\le \dfrac C{|Q|}\dint_1^\ez\dint_Q|f(y)|dy\dfrac {dt}t+\dfrac C{|Q|}\dint_e^\fz\dint_{\{x\in Q:
|f(x)|>t-e\}}|f\chi_Q(y)|dy\dfrac {dt}t\\
&\le \dfrac C{|Q|}\dint_Q M(f\chi_Q)(x)dx+\dfrac C{|Q|}\dint_0^\fz\dint_{\{x\in Q:
|f(x)|>t\}}|f\chi_Q(y)|dy\dfrac {dt}t\\
&\le \dfrac C{|Q|}\dint_Q M(f\chi_Q)(x)dx+\dfrac C{|Q|}\dint_0^\fz|\{x\in Q:\ M(f\chi_Q)(x)>t\}|dt\\
&\le\dfrac C{|Q|}\dint_Q M(f\chi_Q)(x)dx\\
&\le C\dfrac {\Psi(Q)^{\eta/2}}{|Q|}\dint_Q M_{V,\eta/2}(f)(x)dx,
\end{array}$$
since $M(f\chi_Q)(x)\le C\Psi(Q)^{\eta/2}M_{V,\eta/2}(f)(x)$ for
all $x\in Q$. Thus, (4.3) is proved.

Now let us turn to prove
$$M_{\bz,V,\eta}M_{V,\eta+1}f(x)\le CM_{L\log L,\bz,
V,\eta+1}f(x).\eqno(4.4)$$ We first claim that
$$M_{\bz,V,\eta}f(x)\le CM_{L\log L,\bz,
V,\eta+1}f(x).\eqno(4.5)$$ Indeed, it is sufficient to show that
for any cube $Q\ni x$ such that
$$\Psi(Q)^{-\eta}(\Psi(Q)|Q|)^{\bz/n-1}\dint_Q|f(y)|dy\le
\Psi(Q)^{-\eta-1}(\Psi(Q)|Q|)^{\bz/n}\|f\|_{L\log L,Q}.$$ That is,
$|Q|^{-1}\int_Q|f(y)|dy\le \|f\|_{L\log L,Q}$. This is clear from
the definition of the mean Luxemburg norm.

Now we prove (4.4). For any fixed $x\in\rz$ and any fixed cube
$Q\ni x$, write $f=f_1+f_2$, where $f_1=f\chi_{3Q}$. Thus,
$$\begin{array}{cl}
\Psi(Q)^{-\eta}(\Psi(Q)|Q|)^{\bz/n-1}\dint_Q|M_{V,\eta+1}f(y)|dy
&\le\Psi(Q)^{-\eta}(\Psi(Q)|Q|)^{\bz/n-1}\dint_Q|M_{V,\eta+1}f_1(y)|dy\\
&+\Psi(Q)^{-\eta}(\Psi(Q)|Q|)^{\bz/n-1}\dint_Q|M_{V,\eta+1}f_2(y)|dy\\
&:=I+II.
\end{array}$$
For $I$, we know that for all $g$ with $\supp \ g\subset Q$ (see
\cite{p})
$$
\dfrac 1{|Q|}\dint_QMg(y)dy\le C\|g\|_{L\log L,Q}.$$ From this and
note that $M_{V,\eta+1}f(x)\le Mf(x)$, we get
$$\begin{array}{cl}
I&\le
C\Psi(3Q)^{-\eta}(\Psi(3Q)|3Q|)^{\bz/n-1}\dint_{3Q}|Mf_1(y)|dy\\
&\le C\Psi(3Q)^{-\eta-1}(\Psi(Q)|3Q|)^{\bz/n}\|f\|_{L\log
L,3Q}\\
&\le CM_{L\log L,\bz,V,\eta+1}.\end{array}$$
 Next let us estimate $II$. It is easy to see that for all $y,\
 z\in Q$, we have
 $$(\Psi(Q)|Q|)^{\bz/n}M_{V,\eta+1}f_2(y)\le
 CM_{\bz,V,\eta}f(z).\eqno(4.6)$$
In fact,  for any cube $Q'\ni y$ and $Q'\bigcap(\rz\setminus
3Q)\not=\O$, noticing that $z\in Q\subset 3Q'$, we have
$$\begin{array}{cl}
\dfrac
{(\Psi(Q)|Q|)^{\bz/n}}{\Psi(Q')^{\eta+1}|Q'|}\dint_{Q'}|f_2(y)|dy
&\le
C\dfrac {(\Psi(3Q')|3Q'|)^{\bz/n}}{\Psi(3Q')^{\eta+1}|3Q'|}\dint_{3Q'}|f_2(y)|dy\\
&\le CM_{\bz,V,\eta}f(z).\end{array}$$ Hence, (4.6) holds.

Using (4.5) and (4.6) and note that $\Psi(Q)\ge 1$, we obtain
$$\begin{array}{cl}
II&\le C(\Psi(Q)|Q|)^{\bz/n-1}\dint_Q|M_{V,\eta+1}f_2(y)|dy\\
&\le C|Q|^{-1}\dint_Q\dinf_{z\in Q}M_{\bz,V,\eta}f(z)dz\\
&\le CM_{\bz,V,\eta}f(x)\le CM_{L\log
L,\bz,V,\eta+1}f(x).\end{array}$$ Thus, (4.4) is proved.\hfill$\Box$

\begin{lem}\label{l4.2.}\hspace{-0.1cm}{\rm\bf 4.2.}\quad
Let $0\le\bz<n$, $0< \eta<\fz$ and $M_{V,\eta}f, M_{\bz,V,\eta}f$
be locally integrable. Then there is a constant $C>0$ independent
of $f$ and $x$ such that
$$M_{V,\eta+1}M_{\bz,V,\eta}f(x)\le C
M_{\bz,V,\eta}M_{V,\eta}f(x).$$
\end{lem}
{\it Proof}.\quad For any fixed $x\in\rz$ and any fixed cube $Q\ni x$,
write $f=f_1+f_2$, where $f_1=f\chi_{3Q}$. Thus,
$$\begin{array}{cl}
\Psi(Q)^{-\eta-1}|Q|^{-1}\dint_Q|M_{\bz,V,\eta}f(y)|dy
&\le\Psi(Q)^{-\eta-1}|Q|^{-1}\dint_Q|M_{\bz,V,\eta}f_1(y)|dy\\
&\quad+\Psi(Q)^{-\eta-1}|Q|^{-1}\dint_Q|M_{\bz, V,\eta}f_2(y)|dy\\
&:=I+II.
\end{array}$$
Similar to the estimates of $II$ in Lemma 4.1, we have
$$II\le CM_{\bz, V,\eta}f(x)\le CM_{\bz,V,\eta}M_{V,\eta}f(x).$$
From page 881 in \cite{dlz}, we know that for any $g$ with $\supp\
g\subset Q$,
$$\dint_QM_{\bz}f(y)dy\le C|Q|^{\bz/n}\dint_Q g(y)dy,\eqno(4.7)$$
where $M_\bz f(x)=\sup_{x\in Q}\frac 1{|Q|^{1-\bz/n}}\int_Q|f(y)dy$.
Obviously, $M_{\bz,V,\eta}f(x)\le M_{\bz}f(x)$.

Hence, by (4.7), we obtain
$$\begin{array}{cl}
I&\le
C\Psi(3Q)^{-\eta-1}|3Q|^{-1}\dint_{3Q}|M_{\bz,V,\eta}f_1(y)|dy\\
&\le
C\Psi(3Q)^{-\eta-1}|3Q|^{-1}\dint_{3Q}|M_{\bz}f_1(y)|dy\\
&\le
C\Psi(3Q)^{-\eta-1}|3Q|^{\bz/n-1}\dint_{3Q}|f_1(y)|dy\\
&\le
C\Psi(3Q)^{-\eta-1}|3Q|^{\bz/n-1}\dint_{3Q}|M_{V,\eta}f(y)|dy\\
&\le CM_{\bz,V,\eta}M_{V,\eta}f(x).
\end{array}$$
This completes the proof.\hfill$\Box$

\begin{defn}\label{d4.1.}\hspace{-0.1cm}{\rm\bf 4.1.}\quad
Given an increasing function $\vz$, define the function $h_\vz$ by
$$h_\vz(s)=\dsup_{t>0}\dfrac{\vz(st)}{\vz(t)},\quad 0\le s<\fz.$$
\end{defn}

\begin{lem}\label{l4.3.}\hspace{-0.1cm}{\rm\bf 4.3.}\quad
Let $0<\bz<n$, $1\le\eta<\fz$,  $\wz\in A_1^\rho$ and
$B(t)=t\log(e+t)$,. Then there exists a constant $C>0$ such that
for all $t>0$
$$\Phi\l(\wz(\{x\in\rz:\ M_{B,\bz,V,\eta}f(x)>t\})\r)\le
C\dint_\rz\Phi\l(\dfrac {|f(x)|}t\r)h_\Phi(\wz(x))dx,$$ where
$\Phi(s)=s/h_B(s^{\bz/n})$ if $s>0$, otherwise is zero.
\end{lem}
From page 4 in \cite{cf}, we know that
$$\Phi(t)\approx\dfrac {t^{\frac \bz n-1}}{\log(e+t^{\frac
\bz n})}.$$ The the function $\Phi$ is invertible with
$$\Phi^{-1}(t)\approx [t\log(e+t^{\frac \bz n}])^{n/(n-\bz)},$$
and
$$h_\Phi(t)\le C\Theta (t)=Ct^{1-\frac \bz n}\log(e+t^{-\frac
\bz n}).$$
{\it Proof}.\quad  Define $E_t=\{x\in\rz:\
M_{B,\bz,V,\eta}f(x)>t\}$. For each $x\in E_t$, there exists  a
cube $Q_x\ni x$ such that
$$\Psi(Q_x)^{-\eta}(\Psi(Q_x)|Q_x|)^{\bz/n}\|f\|_{B,Q_x}> t.$$
The collection $\{Q_x\}_{x\in E_t}$ covers $E_t$. Similar to the
proof of  Lemma 3.14 in \cite{c}, we can show that there exists a
constant $L> 0$ and a collection of disjoint dyadic cubes
$\{P_j\}$ such that $E_t\subset \bigcup 3P_j$, and such that
$$\Psi(P_j)^{-\eta}(\Psi(P_j)|P_j|)^{\bz/n}\|f\|_{B,P_j}>L t.$$
 By the properties of the Luxemburg norm on Orlicz spaces and by
definition 4.2, we have
$$\begin{array}{cl}
1&\le \dfrac 1{|P_j|}\dint_{P_j}B\l(\dfrac
{(\Psi(P_j)|P_j|)^{\bz/n}|f(x)|}{Lt\Psi(P_j)^{\eta}}\r)dx\\
&\le \dfrac {Ch_B((\Psi(3P_j)|3P_j|)^{\bz/n})}{|3P_j|}
\dint_{P_j}B\l(\frac {|f(x)|}{t\Psi(3P_j)^\eta}\r)dx\\
&\le \dfrac
{Ch_B((\Psi(3P_j)|3P_j|)^{\bz/n})}{\Psi(3P_j)^\eta|3P_j|}
\dint_{P_j}B\l(\frac {|f(x)|}{t}\r)dx\\
&\le \dfrac {Ch_B((\Psi(3P_j)|3P_j|)^{\bz/n})}{\Psi(3P_j)|3P_j|}
\dint_{P_j}B\l(\frac {|f(x)|}{t}\r)dx\\
&= \dfrac {C}{\Phi(\Psi(3P_j)|3P_j|)} \dint_{P_j}B\l(\frac
{|f(x)|}{t}\r)dx.
\end{array}$$
The growth conditions assumed on B imply  that
$$\Phi(\wz(E_t))\le\Phi(\dsum_j(\wz(3P_j))\le
\dsum_j\Phi(\wz(3P_j)).$$
 Hence, if we combine the two inequalities
above and apply Definition 4.2, we get
$$\begin{array}{cl}
\Phi(\wz(E_t))&\le C\dsum_j \dfrac
{\Phi(\wz(3P_j))}{\Phi(\Psi(3P_j)|3P_j|)} \dint_{P_j}B\l(\frac
{|f(x)|}{t}\r)dx\\
&\le C\dsum_j h_\Phi\l(\dfrac {\wz(3P_j)}{\Psi(3P_j)|3P_j|}\r)
\dint_{P_j}B\l(\frac{|f(x)|}{t}\r)dx\\
&\le C\dsum_j\dint_{P_j}B\l(\frac{|f(x)|}{t}\r)h_\Phi(\wz(x))dx\\
&\le C\dint_\rz B\l(\frac{|f(x)|}{t}\r)h_\Phi(\wz(x))dx.
\end{array}$$
Thus, Lemma 4.3 is proved.\hfill$\Box$

\begin{lem}\label{l4.4.}\hspace{-0.1cm}{\rm\bf 4.4.}\quad
Let $b\in BMO(\rho)$, $0<\bz<n$ and $(l_0+1)\le\eta<\fz$. Let
$0<2\dz<\ez<1$, then
$$M^\sharp_{\dz,\eta}([b,I_\bz]f)(x)\le C \|b\|_{BMO(\rho)}(M^\triangle_{\ez, \eta}(I_\bz f)(x)+M_{L\log L,\bz,V,\eta}(f)(x)), \ {\rm a.e\ }\ x\in\rz,\eqno(4.8)$$
holds for any $f\in C_0^\fz(\rz)$.
\end{lem}
{\it Proof}.\quad Observe that for any constant $\lz$
$$[b,I_\bz]f(x)=(b(x)-\lz)I_\bz f(x)-I_\bz((b-\lz)f)(x).$$
As above we fix $x\in\rz$ and let $x\in Q=Q(x_0,r)$(dyadic cube).
Decompose $f=f_1+f_2$, where $f_1=f\chi_{\bar B}$, where $\bar
Q=Q(x,8\sqrt{n}r)$. Let $\lz$ be a constant and $C_Q$ a constant to be
fixed along the proof.

To prove (4.8), we consider two cases about $r$, that is, $r<
\rho(x_0)$ and $r\ge \rho(x_0)$.

Case 1. when  $r< \rho(x_0)$. Since $0<\dz<1$,  we then have
$$\begin{array}{cl}
\l(\dfrac
1{|Q|}\dint_Q|\r.&\l.|[b,I_\bz]f(y)|^\dz-|C_Q|^\dz|\,dy\r)^{1/\dz}\\&
\le \l(\dfrac
1{|Q|}\dint_Q|[b,I_\bz]f(y)-C_Q|^\dz\,dy\r)^{1/\dz}\\
& \le \l(\dfrac
1{|Q|}\dint_Q|(b(y)-\lz)I_\bz f(y)-I_\bz((b-\lz)f)(y)-C_Q|^\dz\,dy\r)^{1/\dz}\\
&\le C\l(\dfrac
1{|Q|}\dint_Q|(b(y)-\lz)I_\bz f(y)|^\dz\,dy\r)^{1/\dz}\\
&\qquad+C\l(\dfrac
1{|Q|}\dint_Q|I_\bz((b-\lz)f_1)(y)|^\dz\,dy\r)^{1/\dz}\\
&\qquad+C\l(\dfrac
1{|Q|}\dint_Q|I_\bz((b-\lz)f_2)(y)-C_Q|^\dz\,dy\r)^{1/\dz}\\
&:=I+II+III.
\end{array}$$
 To deal with $I$, we first fix
$\lz=b_{\bar Q},$ the average of $b$ on $\bar Q$. Then for any
$1<\gz<\ez/\dz$,  note that $m_V(x)\sim m_V(x_0)$ for any $x\in \bar
Q$ and $\Psi(\bar Q)\sim 1$, by  Proposition 4.1, we then obtain
$$\begin{array}{cl}
I &\le C \l(\dfrac 1{|\bar Q|}\dint_{\bar Q}|b(y)-b_{\bar Q}|^{\dz
\gz'}\,dy\r)^{\gz'/\dz}\l(\dfrac
1{|Q|}\dint_Q|I_\bz f(y)|^{\dz \gz}\,dy\r)^{\dz \gz}\\
&\le C \|b\|_{BMO(\rho)}M^\triangle_{\ez,\eta}(I_\bz f)(x),
\end{array}\eqno(4.9)$$ where $1/\gz'+1/\gz=1$.

For II, we recall that $I_\bz$ is weak type $(1,n/(n-\bz))$. Note
that $m_V(x)\sim m_V(x_0)$ for any $x\in \bar Q$ and $\Psi(\bar
Q)\sim 1$, by Kolmogorov's inequality and Proposition 4.3, we then have
$$\begin{array}{cl}
II&\le \dfrac C{|Q|^{1-\bz/n}}\|I_\bz(b-b_{\bar Q}) f_1\|_{L^{\frac n{n-\bz},\fz}}\\
&\le \dfrac C{|\bar Q|^{1-\bz/n}}\dint_{\bar Q}|(b-b_{\bar Q})f(y)|\,dy\\
&\le
CM_{L\log L,\bz,V,\eta}f(x).
\end{array}\eqno(4.10)$$
Finally, for III we first  fix the value of $C_Q$ by taking
$C_Q=(I_\bz((b-b_{\bar Q})f_2))_Q$, the average of
$I_\bz((b-b_{\bar Q})f_2)$ on $B$. Let
$b_{Q_k}=b_{Q(x_0,2^{k+1}r)}$. Then, by Lemmas 2.1 and 2.2, we
have
$$\begin{array}{cl}
II&\le \dfrac C{|Q|}\dint_Q|I_\bz((b-b_{\bar Q})f_2)(y)-(I_\bz((b-b_{\bar Q})f_2))_Q|\,dy\\
&\le \dfrac C{|Q|^2}\dint_Q\dint_Q\dint_{\rz\setminus \bar
Q}|k_\bz(y,\wz)-k_\bz(z,\wz)||(b(\wz)-b_{\bar Q})f(\wz)|d\wz dzdy\\
&\le \dfrac C{|Q|^2}\dint_Q\dint_Q\dint_{|x_0-\wz|>2r
}|k_\bz(y,\wz)-k_\bz(z,\wz)||(b(\wz)-b_{\bar Q})f(\wz)|d\wz dzdy\\
&\le \dfrac
C{|Q|^2}\dint_Q\dint_Q\dsum_{k=2}^\fz\dint_{2^kr\le|x_0-\wz|<2^{k+1}r
}|k_\bz(y,\wz)-k_\bz(z,\wz)||(b(\wz)-b_{\bar Q})f(\wz)|d\wz dzdy\\
&\le C_l\dsum_{k=1}^\fz\dfrac {2^{-
 k\dz_0}}{(1+2^krm_V(x_0))^l|Q_k|^{1-\bz/n}}\dint_{Q_k}|b(\wz)-b_{\bar Q}|f(\wz)|d\wz\\
&\le C_l\dsum_{k=1}^\fz\dfrac
{2^{-k\dz_0}(1+2^krm_V(x_0))^{(l_0+1)\tz(\eta+2)}}{(1+2^krm_V(x_0))^l}\\
&\qquad\qquad\qquad\times\dfrac
1{(1+2^krm_V(Q_k))^{(\eta+1)\az}|Q_k|^{1-\bz/n}}\dint_{Q_k}
|b(\wz)-b_{Q_k}||f(\wz)|d\wz\\
&+ C_l\dsum_{k=1}^\fz\dfrac {2^{-\tz
k}(1+2^krm_V(x_0))^{(l_0+2)\tz(\eta+1)}}{(1+2^krm_V(x_0))^l}\\
&\qquad\qquad\qquad\times\dfrac
1{(1+2^krm_V(Q_k))^{\tz(\eta+2)}|Q_k|^{1-\bz/n}}|b(\bar Q)-b_{Q_k}|\dint_{Q_k}|f(\wz)|d\wz\\
 &\le C_l\dsum_{k=1}^\fz\dfrac
{2^{-
k\dz_0}(1+2^krm_V(x_0))^{(l_0+1)\tz(\eta+2)}}{(1+2^krm_V(x_0))^l}\\
&\qquad\qquad\qquad\times\dfrac
1{(1+2^krm_V(Q_k))^{\tz(\eta+1)}|Q_k|^{1-\bz/n}}\dint_{Q_k}
|b(\wz)-b_{Q_k}||f(\wz)|d\wz\\
&+ C_l\dsum_{k=1}^\fz\dfrac {2^{-\dz_0
k}(1+2^krm_V(x_0))^{(l_0+1)\tz(\eta+2)}}{(1+2^krm_V(x_0))^l}\\
&\qquad\qquad\qquad\times\dfrac
1{(1+2^krm_V(Q_k))^{\tz(\eta+1)}(\Psi(Q_k)|Q_k|)^{1-\bz/n}}|b_{\bar Q}-b_{Q_k}|\dint_{Q_k}|f(\wz)|d\wz\\
&\le C_l\dsum_{k=1}^\fz  2^{-\dz_0 k}\|b\|_{BMO(\rho)}M_{L\log
L,\bz,V,\eta}(f)(x)
+ C_l\|b\|_{BMO(\rho)}M_{\bz,V,\eta}(f)(x)\dsum_{k=1}^\fz k2^{- k\dz_0}\\
&\le C_l\|b\|_{BMO(\rho)}M_{L\log L,\bz,V,\eta}(f)(x),
\end{array}\eqno(4.11)$$
where $l=(l_0+1)(\eta+2)\tz$ and in last inequality we have  used
the fact that  $$|b_{\bar Q}-b_{Q_k}|\le Ck\Psi(Q_k)\|b\|_{BMO(\rho)},$$ and $$ M_{\bz,V,\eta}(f)(x)\le M_{L\log L,\bz,V,\eta}(f)(x).$$
Case 2. When  $r\ge \rho(x_0)$.  Since $0<2\dz< \ez<1$, so
 $a=\eta/\dz$ and $\ez/\dz>2$, then
$$\begin{array}{cl}
\dfrac 1{\Psi(Q)^a}\l(\dfrac
1{|Q|}\dint_Q\r.&\l.|[b,I_\bz]f(y)|^\dz\,dy\r)^{1/\dz}\\
 & =
\dfrac 1{\Psi(Q)^a}\l(\dfrac
1{|Q|}\dint_Q|(b(y)-\lz)I_\bz f(y)-I_\bz((b-\lz)f)(y)|^\dz\,dy\r)^{1/\dz}\\
&\le C\dfrac 1{\Psi(Q)^a}\l(\dfrac
1{|Q|}\dint_Q|(b(y)-\lz)I_\bz f(y)|^\dz\,dy\r)^{1/\dz}\\
&\qquad+C\dfrac 1{\Psi(Q)^a}\l(\dfrac
1{|Q|}\dint_Q|I_\bz((b-\lz)f_1)(y)|^\dz\,dy\r)^{1/\dz}\\
&\qquad+C\dfrac 1{\Psi(Q)^a}\l(\dfrac
1{|Q|}\dint_Q|I_\bz((b-\lz)f_2)(y)|^\dz\,dy\r)^{1/\dz}\\
&:=I+II+III.
\end{array}$$
 To deal with $I$, we first fix
$\lz=b_{\bar Q},$ the average of $b$ on $\bar Q$. Then for any
$2\le\gz<\ez/\dz$, note that $l_0+1\le\eta$, by Proposition 4.1,  we then have
$$\begin{array}{cl}
I &\le C\dfrac 1{\Psi_{\tz'}(Q)} \l(\dfrac 1{|\bar Q|}\dint_{\bar
Q}|b(y)-b_{\bar Q}|^{\dz \gz'}\,dy\r)^{1/(r'\dz)}\\
&\qquad\times\dfrac
{\Psi_{\tz'}(Q)} {\Psi(Q)^{a-\eta/(2\dz)}}\l(\dfrac
1{\Psi(Q)^\eta|Q|}\dint_Q|I_\bz f(y)|^{\dz \gz}\,dy\r)^{1/(\dz \gz)}\\
&\le C \|b\|_{BMO(\rho)}M^\triangle_{\ez,\eta}(I_\bz f)(x),
\end{array}\eqno(4.12)$$ where $1/\gz'+1/\gz=1$.

For II, we recall that $I_\bz$ is weak type $(1,n/(n-\bz))$. By
Kolmogorov's inequality and Proposition 4.3, we then have
$$\begin{array}{cl}
II&\le \dfrac C{\Psi(Q)^a} \dfrac 1{|Q|^{1-\bz/n}}\|I_\bz(b-b_{\bar Q}) f_1\|_{L^{\frac n{n-\bz},\fz}}\\
&\le \dfrac C{\Psi(Q)^{a}}\dfrac 1{|\bar Q|^{1-\bz/n}}\dint_{\bar Q}|(b-b_{\bar Q})f(y)|\,dy\\
&\le CM_{L\log L,\bz,V,\eta}f(x).
\end{array}\eqno(4.13)$$
Finally, for III, let
$b_{Q_k}=b_{Q(x_0,2^{k+1}r)}$. Then, by (1.2) and  Lemma 2.1, we
get that
$$\begin{array}{cl}
III&\le \dfrac C{|Q|}\dint_Q|I_\bz((b-b_{\bar Q})f_2)(y)|\,dy\\
&\le \dfrac C{|Q|}\dint_Q\dint_{\rz\setminus \bar
Q}|k_\bz(y,\wz)||(b(\wz)-b_{\bar Q})f(\wz)|d\wz dy\\
&\le \dfrac C{|Q|}\dint_Q\dint_{|x_0-\wz|>2r
}|k_\bz(y,\wz)||(b(\wz)-b_{\bar Q})f(\wz)|d\wz dy\\
&\le \dfrac
C{|Q|}\dint_Q\dsum_{k=2}^\fz\dint_{2^kr\le|x_0-\wz|<2^{k+1}r
}|k_\bz(y,\wz)||(b(\wz)-b_{\bar Q})f(\wz)|d\wz dy\\
\end{array}$$
$$\begin{array}{cl}
&\le C_l\dsum_{k=1}^\fz\dfrac 1{(1+2^krm_V(x_0))^l|Q_k|^{1-\bz/n}}\dint_{Q_k}|b(\wz)-b_{\bar Q}|f(\wz)|d\wz\\
&\le C_l\dsum_{k=1}^\fz\dfrac
{(1+2^krm_V(x_0))^{(l_0+1)(\eta+2)\tz}}{(1+2^krm_V(x_0))^l}\\
&\qquad\qquad\qquad\times\dfrac
1{(1+2^krm_V(Q_k))^{(\eta+2)\tz}|Q_k|^{1-\bz/n}}\dint_{Q_k}
|b(\wz)-b_{Q_k}||f(\wz)|d\wz\\
&+ C_l\dsum_{k=1}^\fz\dfrac {(1+2^krm_V(x_0))^{(l_0+1)(\eta+2)\tz}}{(1+2^krm_V(x_0))^l}\\
&\qquad\qquad\qquad\times\dfrac
1{(1+2^krm_V(Q_k))^{(\eta+2)\tz}|Q_k|^{1-\bz/n}}|b_{\bar Q}-b_{Q_k}|\dint_{Q_k}|f(\wz)|d\wz\\
 &\le C_l\dsum_{k=1}^\fz\dfrac
{2^{-\dz_0
k}(1+2^krm_V(x_0))^{(l_0+1)(\eta+2)\tz}}{(1+2^krm_V(x_0))^l}\\
&\qquad\qquad\qquad\times\dfrac
1{(1+2^krm_V(Q_k))^{(\eta+2)\tz}|Q_k|^{1-\bz/n}}\dint_{Q_k}
|b(\wz)-b_{Q_k}||f(\wz)|d\wz\\
&+ C_l\dsum_{k=1}^\fz\dfrac {(1+2^krm_V(x_0))^{(l_0+1)(\eta+2)\tz}}{(1+2^krm_V(x_0))^l}\\
&\qquad\qquad\qquad\times\dfrac
1{(1+2^krm_V(Q_k))^{(\eta+1)\tz}(\Psi(Q_k)|Q_k|)^{1-\bz/n}}|b(\bar Q)-b_{Q_k}|\dint_{Q_k}|f(\wz)|d\wz\\
&\le C_l\dsum_{k=1}^\fz  2^{- k\dz_0}\|b\|_{BMO(\rho)}M_{L\log
L,\bz,V,\eta}(f)(x)
+ C_l\|b\|_{BMO(\rho)}M_{\bz,V,\eta}(f)(x)\dsum_{k=1}^\fz k2^{- k\dz_0}\\
&\le C_l\|b\|_{BMO(\rho)}M_{L\log L,\bz,V,\eta}(f)(x),
\end{array}\eqno(4.14)$$  where $l=(l_0+1)(\eta+2)\tz+1$.

 From  (4.9)--(4.14), we get (4.8). Hence the proof is
 finished.\hfill$\Box$

\begin{lem}\label{l4.5.}\hspace{-0.1cm}{\rm\bf 4.5.}\quad
Let $b\in BMO(\rho)$, $0<\bz<n$ and $0<\eta<\fz$. Let $\wz\in A_1^\rho$
and $\Lambda(t)=\Phi(\Phi(t))$. Then there exists a positive
constant $C$ such that for any smooth function $f$ with compact
support
$$\begin{array}{cl}
\dsup_{t>0}&\dfrac 1{\Lambda(1/t)}\wz(\{x\in\rz:\
|[I_\bz,b]f(x)>t\})^{1-\bz/n}\\
&\le C\Lambda(\|b\|_{BMO(\rho)}) \dsup_{t>0}\dfrac
1{\Lambda(1/t)}\wz(\{x\in\rz:\
M_{B,\bz,V,\eta}f(x)>t\})^{1-\bz/n}.\end{array}$$
\end{lem}
Using Theorems 2.1 and Lemmas 4.1, 4.2 and 4.4, adapting the same
arguments in pages 884-886 of \cite{dlz}( or pages 21-23 of
\cite{cf1}), we can prove Lemma 4.5. Here we omit the details.

\bigskip
{\it Proof of Theorem 4.1}.\quad By Lemmas 4.1, 4.4 and Theorems
 3.1, 3.2 and  Corollary 2.1, we have
$$\begin{array}{cl}
\|[I_\bz,b]\|_{L^q(\wz^q)}&\le C\|M_{\ez, V,\eta}(I_\bz
f)\|_{L^q(\wz^q)}+\|M_{L\log L,\bz,V,\eta}(f)\|_{L^q(\wz^q)}\\
&\le C\|I_\bz
f\|_{L^q(\wz^q)}+\|M_{\bz,V,\eta/2} M_{V,\eta/2}f\|_{L^q(\wz^q)}\\
&\le C\|f\|_{L^p(\wz^p)},
\end{array}$$
if taking $\eta=2(p'+q')(l_0+1)$.\hfill$\Box$

\bigskip
{\it Proof of Theorem 4.2}.\quad  By homogeneity, we need only to
show that
$$\wz(\{x\in\rz: \ |[I_\bz,b]f(x)|>1\})\le C\Phi^{-1}\l(\dint_\rz\Phi(
|f(x)|)h_\Phi(\wz(x))dx\r).$$ Indeed, using Lemmas 4.3 and 4.5,
adapting the same arguments pages 21-23 of \cite{cf1}, we can
obtain
$$\begin{array}{cl}
\wz(\{x\in\rz: &\ |[I_\bz,b]f(x)|>1\})\\
&\le C \dsup_{t>0}\dfrac
1{\Lambda(1/t)}\wz(\{x\in\rz:\
|[I_\bz,b]f(x)|>t\})^{1-\bz/n}\\
&\le C\Lambda(\|b\|_{BMO_V}) \dsup_{t>0}\dfrac
1{\Lambda(1/t)}\wz(\{x\in\rz:\
M_{B,\bz,V,\eta}f(x)>t\})^{1-\bz/n}\\
&\le C\Phi^{-1}\l(\dint_\rz\Phi(
|f(x)|)h_\Phi(\wz(x))dx\r),\end{array}$$ if taking $\eta\ge l_0+1$.

 Thus, the proof  is complete.

We next consider a class Schr\"odinger type operators such as
$\nabla(-\Delta+V)^{-1}\nabla$, $\nabla(-\Delta+V)^{-1/2}$,
$(-\Delta+V)^{-1/2}\nabla$  with $V\in B_n$, $(-\Delta+V)^{i\gz}$  with $\gz\in\rr$ and $V\in B_{n/2}$, and $\nabla^2(-\Delta+V)^{-1}$ with $V$ is a nonnegative polynomial.
\begin{thm}\label{t4.3.}\hspace{-0.1cm}{\rm\bf 4.3.}\quad
Let $T$ be operators above, $b\in BMO(\rho)$,  $1<p<\fz$ and $\wz\in A_p^\rho$. Then there exists a constant
$C_p>0$ such that
$$\|[b,T]f\|_{L^p(\wz)}\le C\|b\|_{BMO(\rho)}\|f\|_{L^p(\wz)}.$$
\end{thm}
{\it Proof}.\quad In order to prove Theorem 4.3, from Corollary 2.1, we need only
to show that for
$0<2\dz<\ez<1$ and $\eta>l_0+1$, such that
$$M^\sharp_{\dz,\eta}([b,T]f)(x)\le C \|b\|_{BMO(\rho)}(M_{\ez,V, \eta}(T f)(x)+M_{L\log L,V,\eta}(f)(x)), \ {\rm a.e\ }\ x\in\rz,\eqno(4.15)$$
holds for any $f(x)\in C_0^\fz(\rz)$.\hfill$\Box$

Adapting the same argument of Lemma 4.4, and using (3.1) and (3.2), we can get (4.15). We omit the details here.

The weighted weak-type  endpoint estimate for the commutator of $T$ is the
 following.
\begin{thm}\label{t4.4.}\hspace{-0.1cm}{\rm\bf 4.4.}\quad Let $b\in
BMO(\rho)$ and $\wz\in A_1^{\rho}$. There exists a constant $C>0$ such that for any $\lz>0$
$$\wz(\{x\in\rz: \ |[b,T]f(x)|>\lz\})\le C \dint_\rz\dfrac
{|f(x)|}\lz\l(1+\log^+\l(\frac {|f(x)|}\lz\r)\r)\wz(x)dx.$$
\end{thm}
The proof is similar to that of Theorem 4.1 (see also \cite{p}), we omit the details here.

Finally, we consider another class $V\in B_q$ for $n/2\le q$ for Riesz transforms associated to Schr\"odinger operators. Let $T_1=(-\triangle+V)^{-1}V,\
T_2= (-\triangle+V)^{-1/2} V^{1/2}$ and $T_3=(-\triangle+V)^{-1/2} \nabla$.

\begin{thm}\label{t4.5.}\hspace{-0.1cm}{\rm\bf 4.5.}\quad
Suppose $V\in B_q$ and $q\ge n/2$. Let $b\in BMO(\rho)$. Then
\begin{enumerate}
\item[(i)] If $q'\le p<\fz$ and $\wz\in A_{p/q'}^\rho$,
$$\| [b,T_1]f\|_{L^p(\wz)}\le C\|b\|_{BMO(\rho)}\|f\|_{L^p(\wz)};$$
\item[(ii)] If $(2q)'\le p<\fz$ and $\wz\in A_{p/(2q)'}^\rho$,
$$\| [b,T_2]f\|_{L^p(\wz)}\le C\|b\|_{BMO(\rho)}\|f\|_{L^p(\wz)};$$
\item[(iii)] If $p_0'\le p<\fz$ and $\wz\in A_{p/p_0'}^\rho$, where $1/p_0=1/q-1/n$ and $n/2\le q<n$,
$$\| [b,T_3]f\|_{L^p(\wz)}\le C\|b\|_{BMO(\rho)}\|f\|_{L^p(\wz)}.$$
\end{enumerate}
\end{thm}
The proof of Theorem 4.5 is similar that of Theorem 3.2, we omit the details here.

Let $T_1^*=V(-\triangle+V)^{-1}, T^*_2=  V^{1/2}(-\triangle+V)^{-1/2}$ and $T^*_3=\nabla(-\triangle+V)^{-1/2} $. By duality we can easily get the following results.
\begin{cor}\label{c4.1.}\hspace{-0.1cm}{\rm\bf 4.1.}\quad
Suppose $V\in B_q$ and $q\ge n/2$.  Let $b\in BMO(\rho)$. Then
\begin{enumerate}
\item[(i)] If $1< p\le q$ and $\wz^{-\frac 1{p-1}}\in A_{p'/q'}^\rho$,
$$\|[b, T^*_1]f\|_{L^p(\wz)}\le C\|b\|_{BMO(\rho)}\|f\|_{L^p(\wz)};$$
\item[(ii)] If $1< p\le 2q$ and $\wz^{-\frac 1{p-1}}\in A_{p'/(2q)'}^\rho$,
$$\| [b,T^*_2]f\|_{L^p(\wz)}\le C\|b\|_{BMO(\rho)}\|f\|_{L^p(\wz)};$$
\item[(iii)] If $1< p\le p_0$ and $\wz^{-\frac 1{p-1}}\in A_{p'/p'_0}^\rho$, where $1/p_0=1/q-1/n$ and $n/2\le q<n$,
$$\| [b,T^*_3]f\|_{L^p(\wz)}\le C\|b\|_{BMO(\rho)}\|f\|_{L^p(\wz)}.$$
\end{enumerate}
\end{cor}
We remark that in fact all results in this section also hold for $BMO_{\tz_1}(\rho)$ and $A_p^{\rho,\tz_2}$ if $\tz_1\not=\tz_2$.

 LMAM, School of Mathematical   Sciences

 Peking University

 Beijing, 100871

 P. R. China

\bigskip

 E-mail address:  tanglin@math.pku.edu.cn

\end{document}